\documentclass[11pt]{amsart}
\usepackage{amssymb}

\usepackage{epsfig}

\theoremstyle{plain}
\newtheorem{thm}{Theorem}[section]
\newtheorem*{thm*}{Theorem}
\newtheorem{cor}[thm]{Corollary}
\newtheorem*{cor*}{Corollary}
\newtheorem*{conj*}{Conjecture}
\newtheorem*{lemma*}{Lemma}
\newtheorem{lemma}[thm]{Lemma}
\newtheorem*{prop*}{Proposition}
\newtheorem{prop}[thm]{Proposition}

\theoremstyle{definition}
\newtheorem{rems}[thm]{Remarks}
\newtheorem*{defn*}{Definition}
\newtheorem*{rems*}{Remarks}
\newtheorem*{proof*}{Proof}
\newtheorem{prel*}{Preliminaries}
\newtheorem{examples*}{Examples}

\newcommand{\C}{\mathbb{C}}
\newcommand{\E}{\widetilde{\mathcal E}}
\newcommand{\F}{\widetilde{\mathcal F}}
\newcommand{\M}{\widetilde{M}}
\newcommand{\W}{\widetilde{W}}
\newcommand{\J}{{\mathbb H}}

\newcommand{\npartial}{\not\!\partial}
\newcommand{\SO}{\operatorname{SO}}
\newcommand{\Ind}{\operatorname{Ind}}

\newcommand{\rank}{\operatorname{rank}}
\newcommand{\tr}{\operatorname{tr}}

\newcommand{\R}{\operatorname{\mathbb R}}
\newcommand{\Z}{\operatorname{\mathbb Z}}

\newcommand{\nc}{\newcommand}
\nc{\nt}{\newtheorem}
\nc{\gf}[2]{\genfrac{}{}{0pt}{}{#1}{#2}}
\nc{\mb}[1]{{\mbox{$ #1 $}}}
\nc{\real}{{\mathbb R}}
\nc{\comp}{{\mathbb C}}
\nc{\ints}{{\mathbb Z}}
\nc{\Ltoo}{\mb{L^2({\mathbf H})}}
\nc{\rtoo}{\mb{{\mathbf R}^2}}
\nc{\slr}{{\mathbf {SL}}(2,\real)}
\nc{\slz}{{\mathbf {SL}}(2,\ints)}
\nc{\su}{{\mathbf {SU}}(1,1)}
\nc{\so}{{\mathbf {SO}}}
\nc{\hyp}{{\mathbb H}}
\nc{\disc}{{\mathbf D}}
\nc{\torus}{{\mathbb T}}

\nc{\ca}{{\mathcal A}}
\nc{\cag}{{{\mathcal A}^\Gamma}}
\nc{\cg}{{\mathcal G}}
\nc{\chh}{{\mathcal H}}
\nc{\ck}{{\mathcal B}}
\nc{\cl}{{\mathcal L}}
\nc{\cm}{{\mathcal M}}
\nc{\cs}{{\mathcal S}}
\nc{\cz}{{\mathcal Z}}

\nc{\sind}{\sigma{\rm -ind}}

\begin{document}

\title[TWISTED HIGHER INDEX THEORY ON GOOD ORBIFOLDS]
{TWISTED INDEX THEORY ON GOOD ORBIFOLDS, II:\\
FRACTIONAL QUANTUM NUMBERS}
\author{Matilde Marcolli}
\address{Department of Mathematics, Massachussetts Institute of Technology,
Cambridge, Mass., USA}
\email{matilde@math.mit.edu}
\author{Varghese Mathai}
\address{Department of Mathematics, University of Adelaide, Adelaide 5005,
Australia}
\email{vmathai@maths.adelaide.edu.au}


\begin{abstract}
This  paper uses techniques in noncommutative geometry
as developed by Alain Connes \cite{Co2}, in order
to study the twisted higher index theory
of elliptic operators
on orbifold covering spaces of compact good orbifolds,
which are invariant under a projective
action of the orbifold fundamental group, continuing
our earlier work \cite{MM}.
We also compute the
range of the higher cyclic traces on $K$-theory for
cocompact Fuchsian groups, which is then applied to determine the
range of values of the Connes-Kubo Hall conductance in the
discrete model of the quantum Hall effect on the  hyperbolic
plane, generalizing earlier results in
\cite{Bel+E+S}, \cite{CHMM}.  The new phenomenon that we observe in
our case is that the Connes-Kubo Hall conductance has plateaux
at integral multiples of a {\em fractional valued topological
invariant}, namely the {\em orbifold Euler characteristic}.
Moreover the set of possible fractions has been determined, and
is compared with recently available experimental data.
It is plausible that this might shed some light on the
mathematical mechanism responsible for fractional quantum numbers.
\end{abstract}

\maketitle

\section*{Introduction}
This  paper uses techniques in noncommutative geometry
as developed by Alain Connes \cite{Co2}
in order to prove a twisted higher
index theorem for elliptic operators
on orbifold covering spaces of compact good orbifolds,
which are invariant under a projective
action of the orbifold fundamental group. These
higher indices are basically
the evaluation of pairings of higher traces (which are
cyclic cocycles arising from the orbifold fundamental group
and the multiplier defining the projective action) with
the index of the elliptic operator, considered as an element
in the $K$-theory of some completion of the twisted
group algebra of the orbifold fundamental group.
This paper is the continuation of \cite{MM} and generalizes
the results there.
The main purpose for studying the twisted higher index
theorem on orbifolds is to highlight the fact that
when the orbifold is not smooth, then the twisted higher
index can be a  {\em fraction}. In particular, we determine
the range of the higher cyclic traces on $K$-theory
for general cocompact Fuchsian groups.
We adapt and generalize the discrete model
of the quantum Hall effect of Bellissard
and his collaborators \cite{Bel+E+S} and also \cite{CHMM},
to the case of general cocompact Fuchsian
groups and orbifolds, which can  be viewed equivalently
as the generalization to the equivariant context.
The new phenomenon that we observe in
our case is that the Connes-Kubo Hall conductance has plateaux
at integral multiples of a {\em fractional} valued {\em topological
invariant}, namely the {\em orbifold Euler characteristic}.
The presence of denominators is
caused by the presence of cone points singularities
and by the hyperbolic
geometry on the complement of these cone points.
The negative
curvature of the hyperbolic structure replaces interaction and
simulates, in our single electron model, the presence of Coulomb
interactions. We also have
a geometric term in the Hamiltonian (arising from the cone point
singularities) which accounts partly for the effect of Coulomb
interactions.
This geometric model of interaction is fairly simple, hence the
agreement of our fractions with the experimental values is only
partial. Among the observed fractions, for instance, we can derive
$5/3$,
$4/3$,
$4/5$, $2/3$, and $5/2$ from genus one orbifolds,  and $2/5$,
$1/3$, $4/9$,
$4/7$, $3/5$,
$5/7$, $7/5$ from genus zero orbifolds, see \S 5. However,
fractions like $3/7$ and
$5/9$, seem unobtainable in  this model, even including
higher genus orbifolds. Their explanation probably requires a more
sophisticated term describing the electron interaction. It is not
unreasonable to expect that this term may also be geometric in
nature, but we leave it to future studies.

There are currently several different models which describe
the occurrence of fractional quantum numbers in the quantum Hall
effect. Usually quantum field theoretic techniques are involved. Most
notably, there is a sophisticated Chern-Simons theory
model for the fractional quantum Hall effect developed by
Frohlich and his collaborators, cf.
\cite{Froh}. Also within the quantum field theoretic formalism it can be
noticed that possibly different models are needed in order
to explain the occurrence of different sets of fractions. For example,
the fraction $5/2$ requires by itself a separate model.

After reviewing some preliminary material in section 1,
we establish in section 2 a twisted higher index theorem which
adapts the proofs of the index
theorems of Atiyah \cite{At},  Singer \cite{Si},
Connes and Moscovici \cite{CM}, and Gromov
\cite{Gr2}, \cite{Ma1}, to the case of
good orbifolds, that is, orbifolds whose
orbifold universal cover is a smooth manifold.
This theorem generalizes the twisted index theorem for $0$-traces
of \cite{MM} to the case of higher degree cyclic traces. The result
can be summarized as follows.
Let
$\mathcal{R}$ be the algebra
of rapidly decreasing sequences, i.e. \[
\mathcal{R}=\left\{\left(a_{i}\right)_{i\in\mathbb{N}}\,:\,
\sup_{i\in\mathbb{N}} i^k  \left|a_{i}\right|<\infty
\ \forall\ k\in\mathbb{N}\right\}
\] Let $\Gamma$ be a discrete group and
$\sigma$ be a multiplier on $\Gamma$. Let $\mathbb{C}(\Gamma, \sigma)$
denote the twisted group algebra. We denote the tensor product
$\mathbb{C}(\Gamma, \sigma) \otimes \mathcal R$ by $\mathcal R(\Gamma,
\sigma)$.
Let $\Gamma\to\widetilde{M}\to M$ denote the universal orbifold
cover of a compact good orbifold $M$, so that $\widetilde M$ is a
smooth manifold.
Suppose given a multiplier $\sigma$ on
$\Gamma$ and assume that there is a projective $(\Gamma,\bar\sigma)$-action
on $L^2$ sections of $\Gamma$-invariant vector bundles over $\widetilde M$.
By considering $(\Gamma,\bar\sigma)$-invariant elliptic operators $D$ acting on
$L^2$ sections of these bundles, we will define a
$(\Gamma, \sigma)$-index element in $K$-theory
\[
\Ind_{\sigma} ({D})\in
K_0(\mathcal{R}(\Gamma,\sigma)). \] We will compute the
pairing of
$\Ind_{\sigma}({D})$ with higher traces. More
precisely, given a normalized group cocycle $c\in Z^k(\Gamma,
\mathbb C)$, we define a cyclic cocycle
$\tr_c \in ZC^k(\C(\Gamma, \sigma))$ of dimension $k$ on the twisted group
algebra
$\mathbb C(\Gamma, \sigma)$, which extends
continuously to a $k$-dimensional cyclic cocycle on $\mathcal{R}(\Gamma,
\sigma)$. This induces a map on $K$-theory,
$$
[\tr_c] : K_0 (\mathcal{R}(\Gamma, \sigma) )  \to \C.
$$
A main theorem established in this paper is a cohomological formula for
$$
{\Ind}_{(c, \Gamma ,\sigma)}({D}) =
[\tr_c]\left(\Ind_{\sigma}({D})\right).
$$
Our method consists of applying the Connes-Moscovici
{\em local} higher index theorem to a
family of idempotents constructed from the heat operator on $\widetilde{M}$,
all of which represent the $(\Gamma, \sigma)$-index.

Let $\Gamma$ be a Fuchsian group of signature $(g; \nu_1,\ldots, \nu_n)$,
that is, $\Gamma$ is the orbifold fundamental group of
the 2 dimensional hyperbolic orbifold $\Sigma(g; \nu_1,\ldots,\nu_n)$
of signature $(g; \nu_1,\ldots,\nu_n)$, where $g \in \Z, g\ge 0$
denotes the genus and $2\pi  /\nu_j,\; \nu_j \in \mathbb N$ denotes
the cone
angles at the cone points of the orbifold.
In \cite{MM} we computed
the $K$-theory of the twisted group $C^*$ algebra. Under the
assumption that the Dixmier-Douady invariant of the multiplier
$\sigma$ is trivial, we obtained
$$
K_j(C^*(\Gamma, \sigma)) \cong \left\{\begin{array}{l}
\Z^{2-n +\sum_{j=1}^n \nu_j}\qquad \text{if}\,\, j=0;\\ \\
\mathbb{Z}^{2g} \qquad \qquad\qquad\text{if}\,\, j=1.
\end{array}\right.
$$
Here we use a result of \cite{Ji},
which is a twisted analogue of a result of Jollissant
and which says in particular that, when $\Gamma$ is a cocompact Fuchsian
group, then
the natural inclusion map $j: \mathcal{R}(\Gamma,\sigma) \to
C^*(\Gamma,\sigma)$ induces an isomorphism in $K$-theory
$$
K_\bullet (\mathcal{R}(\Gamma,\sigma) ) \cong K_\bullet (C^*(\Gamma,\sigma))
$$
Using this, together with our twisted higher index theorem for good
orbifolds and some results in \cite{MM}, and under the same assumptions as
before,
we determine, in section 3, the range of the higher trace
on
$K$-theory
$$
[\tr_c]( K_0 (C^*(\Gamma,\sigma))) = \phi\Z,
$$
where $ -\phi = 2(1-g) + (\nu-n)\in {\mathbb Q}$ is the orbifold Euler
characteristic of $\Sigma(g;\nu_1, \ldots, \nu_n)$. Here we have
 $\quad\nu = \sum_{j=1}^n 1/\nu_j$ and
$c$ is the area 2-cocycle on $\Gamma$, i.e. $c$ is the restriction to
$\Gamma$
of the area 2-cocycle on $PSL(2, \R)$.

In section 4  we study
the hyperbolic Connes-Kubo formula for the
Hall conductance in
the discrete model of the Quantum Hall Effect on the hyperbolic plane,
where we consider Cayley graphs of Fuchsian groups which may have
torsion subgroups. This generalizes the results in \cite{CHMM} where
only torsion-free Fuchsian groups were considered. We recall that
the results in \cite{CHMM} generalized to
hyperbolic space the noncommutative geometry
approach to the Euclidean
quantum Hall effect that was pioneered by
Bellissard and collaborators \cite{Bel+E+S}, Connes \cite{Co}
and Xia \cite{Xia}.
We first relate the hyperbolic Connes-Kubo Hall conductance cyclic
2-cocycle and the area cyclic 2-cocycle on the algebra ${\mathcal
R} (\Gamma,
\sigma)$,
and show that they define the same class in cyclic cohomology.
Then we use our theorem on the range of the higher trace on $K$-theory to
determine
the range of values of the Connes-Kubo Hall conductance cocycle in
the Quantum Hall Effect.
The new phenomenon
that we observe in this case is that the Hall conductance
has plateaux at all energy levels belonging to
any gap in the spectrum of the Hamiltonian
(known as the generalized
Harper operator), where it
is now shown to be equal to an integral multiple of
a {\em fractional} valued
topological invariant $\phi$, which is the negative of
the orbifold Euler characteristic of the good orbifold
$\Sigma(g;\nu_1,\ldots,\nu_n)$. If we fix the genus, then
the set of possible denominators is finite by the Hurwitz theorem \cite{Sc},
and has been explicitly determined in the low genus cases
\cite{Bro}. This provides a topological explanation of the appearance of
fractional quantum numbers. In the last section we compare our
results with some observed values.

In section 5, we provide lists of specific examples of
good 2-dimensional orbifolds for which $\phi$ is {\em not} an integer.
First we observe how
the presence of both the hyperbolic structure and  the
cone points is essential in order to have fractional quantum
numbers. In fact, $\phi$ is an
{integer} whenever the hyperbolic orbifold is smooth, i.e. whenever
$1=\nu_1=\ldots =\nu_n$, which is the case considered in
\cite{CHMM}. Similarly, by direct inspection, it is possible
to see that all Euclidean orbifolds also produce only {\em integer}
values of $\phi$. We use the class of orbifolds which are spheres or
tori with cone points,
having a (singular) hyperbolic structure, to represent in our physical
model some of
the fractions observed in the FQHE.
We also list the examples arising from quotients of low genus surfaces
\cite{Bro}, and we discuss some phenomenology on the role of the
orbifold points and of the minimal genus $g'$ of the covering surface.

Summarizing, one key advantage of our model is that the fractions
we get are obtained from an equivariant index theorem
and are thus topological in nature. Consequently, as pointed
out in \cite{Bel+E+S}, the Hall conductance
is seen to be stable under small deformations of the Hamiltonian.
Thus, this model can be easily generalized to systems with
disorder as in \cite{CHM}. This is a necessary step
in order to establish the presence of plateaux \cite{Bel+E+S}. The
main limitation of our model is that there is a small number of
experimental fractions that we do not obtain in our model, and we also
derive other fractions which do not seem to correspond to experimentally
observed values.
To our knowledge, however, this is also a limitation occuring in the
other models available in the literature.

\noindent{\bf Acknowledgments:}
We thank J. Bellissard for his encouragement and  for some useful
comments. The second author thanks A. Carey and K.
Hannabuss for some helpful comments concerning the section 4. The
first author is partially supported by NSF grant DMS-9802480. Research by
the second author is supported by the Australian Research Council.

\section{Preliminaries}

Recall that, if $\mathbb H$ denotes the hyperbolic
plane and $\Gamma$ is a Fuchsian group of signature $(g; \nu_1,\ldots,
\nu_n)$, that is, $\Gamma$ is a discrete cocompact subgroup of
$PSL(2, \R)$ of genus $g$ and with $n$ elliptic elements of order
$ \nu_1,\ldots, \nu_n$ respectively, then
the corresponding compact oriented hyperbolic 2-orbifold of signature
$(g; \nu_1,\ldots, \nu_n)$ is defined as the quotient space
$$
\Sigma(g; \nu_1,\ldots,\nu_n)= \Gamma\backslash {\mathbb H},
$$ where $g$ denotes the genus and $2\pi /\nu_j,\; \nu_j \in
\mathbb N$
denotes the cone angles at the cone points of the orbifold.
A compact oriented 2-dimensional Euclidean orbifold is obtained
in a similar manner, but with $\mathbb H$ replaced by $\R^2$.

All Euclidean and hyperbolic 2-dimensional
orbifolds $\Sigma(g;\nu_1,\ldots,\nu_n)$ are
good, being in fact orbifold covered by a smooth surface
$\Sigma_{g'}$ cf. \cite{Sc}, i.e. there is a finite group $G$ acting on
$\Sigma_{g'}$ with quotient $\Sigma(g;\nu_1,\ldots,\nu_n)$,
where $g' =1+ \frac{\#(G)}{2}(2(g -1) + (n - \nu))$ and where
$\nu = \sum_{j=1}^n 1/\nu_j$.

For fundamental material on orbifolds, see \cite{Sc}, \cite{FuSt} and
\cite{Bro}. See also \cite{MM}, section 1.

Let $M$ be a good, compact orbifold, and $\mathcal E \to M$ be an orbifold
vector bundle over $M$, and $\E \to \widetilde{M}$ be its lift to
the universal orbifold covering space $\Gamma\to \widetilde{M}\to M$,
which is by assumption
a simply-connected smooth manifold. We have a
$(\Gamma, \bar\sigma)$-action (where $\sigma$ is a multiplier on
$\Gamma$ and $\bar\sigma$ denotes its complex conjugate)
on $L^2(\M)$,where we choose $\omega = d\eta$ an exact $2$-form on
$\widetilde{M}$
such that $\omega$ is also $\Gamma$-invariant, although $\eta$ is {\em not}
assumed to be $\Gamma$-invariant, and the Hermitian connection
\[
   \nabla = d + i\eta
\]
on the trivial line bundle over $\widetilde{M}$,
with curvature $\nabla^2=i\omega$. The projective action
is defined as follows:\\

Firstly, observe that since $\widetilde{\omega}$ is $\Gamma$-invariant,
$0=\gamma^*\widetilde{\omega}-\widetilde{\omega}=d(\gamma^*\eta-\eta)\
\forall\gamma \in\Gamma$. So $\gamma^*\eta-\eta$ is a closed 1-form on the
simply connected manifold $\widetilde{M}$, therefore \[
\gamma^*\eta-\eta=d\phi_\gamma\quad\forall\gamma\in\Gamma \]
where $\phi_\gamma$ is a smooth function on $\widetilde{M}$ satisfying in
addition,
\begin{itemize}
\item $\phi_\gamma(x)+\phi_{\gamma'}(\gamma x)-\phi_{\gamma'\gamma}(x)$
is independent of $x\in\widetilde{M}\ \forall \gamma,\gamma'\in\Gamma$;
\item $\phi_\gamma(x_0)=0$ for some $x_0\in\widetilde{M}\ \forall\in\Gamma$.
\end{itemize}
Then $\bar{\sigma}(\gamma,\gamma')=\exp(i\phi_\gamma(\gamma'\cdot x_0))$
defines a multiplier on $\Gamma$ i.e. $\bar{\sigma}:\Gamma\times\Gamma\to U(1)$
satisfies the following identity for all $\gamma, \gamma', \gamma'' \in \Gamma$
\[ \bar\sigma(\gamma,\gamma')\bar\sigma(\gamma,\gamma'\gamma'')=
\bar\sigma(\gamma\gamma',\gamma'') \bar\sigma(\gamma',\gamma'') \]
For $u \in L^2(\widetilde{M},
\E)$, let $
S_\gamma u  =  e^{i\phi_\gamma}u$ and $U_\gamma u =
\gamma^*u$
and $T_\gamma=U_\gamma \circ S_\gamma$ be the composition.
Then $T$ defines a projective $(\Gamma,\bar{\sigma})$-action on $L^2$-spinors,
i.e. \[
T_\gamma T_{\gamma'}=\bar{\sigma}(\gamma,\gamma')T_{\gamma\gamma'}. \]

This defines
a $(\Gamma,\bar{\sigma})$-action, provided that the Dixmier-Douady invariant
$ \delta(\sigma) = 0$, see \cite{MM}.

As in \cite{MM}, we shall consider the {\em twisted group von Neumann
algebra} $W^*(\Gamma,\sigma)$, the commutant of the left
$\bar{\sigma}$-regular representation on $\ell^2(\Gamma)$ and $W^*( \sigma)
$ as the
commutant of the $(\Gamma, \bar\sigma)$-action on
$L^2(\widetilde{M}, \widetilde{\mathcal{S}^\pm\otimes E})$.

We have an identification (see \cite{MM})
$$
W^*( \sigma) \cong
W^*(\Gamma, \sigma) \otimes B(L^2({\mathcal F}, \E|_{{\mathcal F}}) ) $$
where $B(L^2({\mathcal F}, \E|_{{\mathcal F}}) ) $ denotes the  algebra of all
bounded operators
on the Hilbert space $L^2({\mathcal F}, \E|_{{\mathcal F}}) $, and
${\mathcal F}$ is a relatively compact fundamental domain in
$\M$ for the action of $\Gamma$. We have a
semifinite trace
$$
\tr: W^*(\sigma) \to \C
$$
defined as in the untwisted case due  to Atiyah \cite{At},
$$
Q \to \int_{M} \tr(k_Q(x,x))dx
$$
where $k_Q$ denotes the Schwartz kernel of $Q$. Note that this trace is finite
whenever $k_Q$ is continuous in a neighborhood of the diagonal in
${\widetilde M}\times {\widetilde M}$.

We also consider, as in \cite{MM}, the subalgebra
${C^*}( \sigma)$ of $W^*(\sigma)$, whose elements have the additional
property of
some off-diagonal decay, and one also has the identification (cf. \cite{MM})
$$
C^*( \sigma) \cong
C^*(\Gamma, \sigma) \otimes {\mathcal K}(L^2({\mathcal F}, \E|_{{\mathcal
F}}) )
$$

In \cite{MM} we considered the $C^*$ algebra
$$ C^*(M)=C(P)
\rtimes \SO(m), $$
where $P$
is the bundle of oriented frames on the orbifold tangent bundle.
The relevent K-theory is the {\em orbifold $K$-theory}
$$K^0_{orb}(M) \equiv K_0(C^*(M))
= K_0(C(P) \rtimes \SO(m)) \cong K^0_{\SO(m)}(P). $$
In the case when $M$ is a good orbifold, one can show
that the $C^*$ algebras $C^*(M)$ and $C_0(X)\rtimes G$ are
strongly Morita equivalent,
where $X$ is smooth and $G\to X \to M$ is an orbifold cover.
In particular,
$$K^0_{orb}(M) \cong K^0(C_0(X)
\rtimes G) = K^0_G (X).$$

The relevant cohomology is the {\em orbifold cohomology} $ H^j_{orb}(M) =
H^j(X, G)$, for $j=0,1$, which is the delocalized equivariant cohomology for
a finite group action on a smooth manifold \cite{BC}. The
Baum-Connes equivariant Chern character is a homomorphism
$$ ch_G : K^0_G(X)\to H^0(X,G). $$

Let $\underline B\Gamma = \Gamma\backslash \underline E\Gamma$ be the
classifying space of proper actions, as defined in
\cite{BCH}. In our case, the orbifold
$\Sigma(g;\nu_1,\ldots,\nu_n)$,
viewed as the quotient space
$\Gamma\backslash {\mathbb H}$, is
$\underline B\Gamma(g;\nu_1,\ldots,\nu_n)$.
Equivalently, $\underline B\Gamma(g;\nu_1,\ldots,\nu_n)$ can be viewed
as the classifying space of the orbifold fundamental group
$\Gamma(g;\nu_1,\ldots,\nu_n)$.

Let $S\Gamma$ denote the set of all elements of $\Gamma$ which are of finite
order. Then $S\Gamma$ is not empty, since $1\in S\Gamma$. $\Gamma$
acts on $S\Gamma$ by conjugation, and let $F\Gamma$ denote
the associated permutation module over $\C$, i.e.
$$
F\Gamma = \left\{\sum_{\alpha \in S\Gamma} \lambda_\alpha [\alpha]
\, \Big|\,
\lambda_\alpha \in \C \quad \text{and}\quad \lambda_\alpha = 0
\quad \text{except for
a finite number of}\quad\alpha\right\}
$$

Let $C^k(\Gamma, F\Gamma) $ denote the space of all antisymmetric
$F\Gamma$-valued $\Gamma$-maps on $\Gamma^{k+1}$, where $\Gamma$ acts
on  $\Gamma^{k+1}$ via the diagonal action. The coboundary map is
$$
\partial c(g_0, \ldots, g_{k+1}) =
\sum_{i=0}^{k+1}
(-1)^i c(g_0, \ldots,\hat g_i \ldots g_{k+1})
$$
for all $c \in C^k(\Gamma, F\Gamma)$ and
where $\hat g_i $ means that $g_i$ is omitted.
The cohomology
of this complex is the group cohomology of $\Gamma$
with coefficients in $F\Gamma$,
$H^k(\Gamma, F\Gamma) $, cf. \cite{BCH}. They also show that
$H^k(\Gamma, F\Gamma) \cong H^j(\Gamma, \C)\oplus_m
H^k(Z(C_m),\C)$, where
$S\Gamma = \{1, C_m | m = 1, \ldots \}$ and the isomorphism is canonical.

Also, for any Borel measurable $\Gamma$-map
$\mu: \underline E\Gamma \to \Gamma$, there is an induced map
on cochains
$$
\mu^*: C^k(\Gamma, F\Gamma) \to C^k(\underline E\Gamma, \Gamma)
$$
which induces an isomorphism on cohomology,
$\quad \mu^* : H^k(\Gamma, F\Gamma)
\cong H^k(\underline E\Gamma, \Gamma)\quad$
\cite{BCH}. Here $H^j(\underline E\Gamma, \Gamma)$
denotes the $\Z$-graded (delocalised) equivariant cohomology
of $\underline E\Gamma$, which is a refinement of what was
discussed earlier, and which is
defined in \cite{BCH} using sheaves (and cosheaves), but we  will not
recall the definition here.

Let $M$ be a good orbifold with orbifold fundamental group $\Gamma$.
We have seen that the universal orbifold cover $\widetilde M$
is classified by a
continuous map $f: M \to \underline B\Gamma$, or equivalently
by a $\Gamma$-map $f: \widetilde M \to \underline E\Gamma$.
The induced map is $f^* :
H^j_{orb}(\underline B\Gamma, \C)  \equiv
H^k(\underline E\Gamma, \Gamma)
\to  H^k(\widetilde M, \Gamma) \equiv H^k_{orb}( M, \C)$ and
therefore in particular one has $f^*([c]) \in  H^k_{orb}( M, \C)$
for all $[c] \in H^k(\Gamma, \C)$.
This can be expressed on the level of cochains by easily
modifying
the procedure in \cite{CM}, and we refer to
\cite{CM} for further details.

Finally, we add here a brief comment on the assumption used throughout
\cite{MM} on the vanishing of the Dixmier-Douady invariant of the
multiplier $\sigma$. We show here that the condition is indeed necessary,
since we can always find examples where $\delta(\sigma)\neq 0$.
Let $\Gamma$ be the Fuchsian group of signature
$(g;\nu_1,\ldots,\nu_n)$, as before.
Consider the long exact sequence of the change of coefficient groups,
as in \cite{MM},
\begin{equation}
\label{change}
\begin{array}{c}
\cdots H^1(\Gamma,U(1))\stackrel{\delta}\to
H^2(\Gamma,\Z)\stackrel{i_*}{\to}
H^2(\Gamma,\R)\stackrel{e^{2\pi\sqrt{-1}}_*}{\to} \\[2mm]
 H^2(\Gamma,U(1)) \stackrel{\delta}{\to} H^3(\Gamma,\Z) \to
H^3(\Gamma,\R).
\end{array}\end{equation}
The argument of \cite{MM} shows that $H^3(\Gamma,\R)=0$ and
$H^2(\Gamma,\R)=\R$. Moreover, we observe that $H^1(\Gamma,\Z)=
Hom (\Gamma, \Z) \cong \Z^{2g}$, $\; H^1(\Gamma,\R)=
Hom (\Gamma, \R) \cong \R^{2g}$ and
$\; H^1(\Gamma,U(1))= Hom (\Gamma, U(1)) \cong U(1)^{2g}\times_{j=1}^n
\Z_{\nu_j}$. Now
$H^2(\Gamma,\Z)=\Z\oplus_j \Z_{\nu_j}$, see \cite{Patt}, which is
consistent with
the result in \cite{MM} that the group of the orbifold line bundles over the
orbifold $\Gamma\backslash {\J}$ has $1-n + \sum_{j=1}^n {\nu}_j$ generators.
It is also proved in \cite{Patt} that
$H^2(\Gamma, U(1)) = U(1)\times_{j=1}^n \Z_{\nu_j}$. Using
the long exact sequence and the remarks above,
we see that $H^3(\Gamma, \Z) = Tor(H^2(\Gamma, U(1))) = \times_{j=1}^n
\Z_{\nu_j}$. Thus, in the sequence we have $Ker(i_*)=\oplus_j
\Z_{\nu_j}$, $Im(i_*)=\Z=Ker(e^{2\pi\sqrt{-1}}_*)$,
$Im(e^{2\pi\sqrt{-1}}_*)=U(1)$. So we can identify all the classes of
multipliers with trivial Dixmier--Douady invariant with
$U(1)=Ker(\delta)$. Finally, we have
$$ Im(\delta)=H^3(\Gamma,\Z)=H^2(\Gamma,U(1))/Ker(\delta)= \oplus_j
\Z_{\nu_j}.  $$

The calculations of the cohomology of the Fuchsian group $\Gamma =
\Gamma(g; \nu_1, \ldots, \nu_n)$ are summarized in the following table.

{

\begin{tabular}{|c|c|c|c|} \hline
$j $ &   $\qquad\qquad H^j(\Gamma, \Z)\qquad\qquad$ &
$\qquad\qquad H^j(\Gamma, \R)\qquad\qquad$ &  $\qquad\qquad
H^j(\Gamma, {\bf U}(1))\qquad\qquad$\\
\hline \hline
$0$  &     $\Z$  &   $\R$ &  ${\bf U}(1)$\\  \hline
$1$ &   $\Z^{2g}$ &   $\R^{2g}$&   ${\bf U}(1)^{2g}\oplus_j\Z_{\nu_j}$ \\
\hline
$2$ &    $\Z\oplus_j\Z_{\nu_j}$ &   $\R$ &   ${\bf U}(1)\oplus_j\Z_{\nu_j}$
\\ \hline
$3$ &   $\oplus_j\Z_{\nu_j}$ &   $0$ &   $$ \\ \hline
\end{tabular}

}

\section{ Twisted higher index theorem}

In this section, we will define the higher twisted index of an elliptic
operator on a good orbifold, and establish a cohomological formula for
any cyclic trace arising from a group cocycle, and
which is applied to the twisted higher index. We adapt the strategy
and proof in \cite{CM} to our context.

\subsection{Construction of the parametrix and the index map}
Let $M$ be a compact, good orbifold, that is, the universal cover
$\Gamma \to \M \to M$ is a smooth manifold
and we will assume, as before, that there is a $(\Gamma, \bar\sigma)$-action
on $L^2(\M)$ given by $T_\gamma = U_\gamma \circ S_\gamma \, \forall \gamma
\in \Gamma$. Let $\E, \ \F$ be Hermitian vector bundles on $M$ and let
$\E, \ \F$ be the corresponding lifts to $\Gamma$-invariants Hermitian vector
bundles on $\M$. Then there are induced $(\Gamma, \bar\sigma)$-actions
on $L^2(\M, \E)$ and $L^2(\M, \F)$ which are also given by
$T_\gamma = U_\gamma \circ S_\gamma \, \forall \gamma \in \Gamma$.

Now let $D : L^2(\M, \E) \to L^2(\M, \F)$ be a first order
$(\Gamma, \bar\sigma)$-invariant elliptic operator. Let $U\subset \M$
be an open subset that contains the closure of a fundamental domain
for the $\Gamma$-action on $\M$. Let $\psi \in C^\infty_c(\M)$ be a
compactly supported smooth function such that $\text{supp}(\psi) \subset
U$, and
$$
\sum_{\gamma \in \Gamma} \gamma^*\psi = 1.
$$
Let $\phi \in C^\infty_c(\M)$ be a
compactly supported smooth function such that $\phi=1$
on $\text{supp}(\psi)$.

Since $D$ is elliptic, we can construct a parametrix $J$ for it on the open
set $U$ by standard methods,
$$
JDu = u - Hu \qquad \forall u \in C^\infty_c(U, \E|_U)
$$
where $H$ has a smooth Schwartz kernel.
Define the pseudodifferential operator $Q$ as
\begin{equation}
\label{eqno(1)}
Q = \sum_{\gamma \in \Gamma} T_\gamma \phi J \psi T_\gamma^*
\end{equation}
We compute,
\begin{equation}\label{eqno(2)}
QD w = \sum_{\gamma \in \Gamma} T_\gamma \phi J \psi D T_\gamma^* w \qquad
\forall w \in C^\infty_c(\M, \E),
\end{equation}
since $T_\gamma D = D T_\gamma \quad \forall\gamma\in \Gamma$. Since $D$ is
a first order operator, one has
$$
D(\psi w ) = \psi Dw + (D\psi) w
$$
so that (\ref{eqno(2)}) becomes
$$
=\sum_{\gamma \in \Gamma} T_\gamma \phi J D \psi  T_\gamma^* w
- \sum_{\gamma \in \Gamma} T_\gamma \phi J (D\psi)  T_\gamma^* w.
$$
Using (\ref{eqno(1)}), the expression above becomes
$$
= \sum_{\gamma \in \Gamma} T_\gamma \psi  T_\gamma^* w
- \sum_{\gamma \in \Gamma} T_\gamma \phi H \psi  T_\gamma^* w
- \sum_{\gamma \in \Gamma} T_\gamma \phi J (D \psi)  T_\gamma^* w.
$$
Therefore (\ref{eqno(2)}) becomes
$$
QD = I - R_0
$$
where
$$
R_0 = \sum_{\gamma \in \Gamma} T_\gamma \left(\phi H \psi
+  J (D \psi) \right) T_\gamma^*
$$
has a smooth Schwartz kernel.
It is clear from the definition that one has $T_\gamma Q = Q T_\gamma$
and  $T_\gamma R_0 = R_0 T_\gamma \quad \forall\gamma\in \Gamma$.
Define
$$
R_1 = ^tR_0 + D R_0 ^tQ - DQ(^t R_0).
$$
Then  $T_\gamma R_1 = R_1 T_\gamma \quad \forall\gamma\in \Gamma$,
$R_1$ has a smooth Schwartz kernel and satisfies
$$
DQ = I - R_1.
$$
Summarizing, we have obtained the following

\begin{prop}
Let $M$ be a compact, good orbifold and $\Gamma \to \M \to M$
be the universal orbifold
covering space. Let $\mathcal E, \ \mathcal F$ be Hermitian vector bundles on
$M$ and let
$\E, \ \F$ be the corresponding lifts to $\Gamma$-invariants Hermitian vector
bundles on $\M$. We will assume as before that there is a $(\Gamma,
\bar\sigma)$-action
on $L^2(\M)$ given by $T_\gamma = U_\gamma \circ S_\gamma \, \forall \gamma
\in \Gamma$, and induced $(\Gamma, \bar\sigma)$-actions
on $L^2(\M, \E)$ and $L^2(\M, \F)$ which are also given by
$T_\gamma = U_\gamma \circ S_\gamma \, \forall \gamma \in \Gamma$.

Now let $D : L^2(\M, \E) \to L^2(\M, \F)$ be a first order
$(\Gamma, \bar\sigma)$-invariant elliptic operator. Then there
is an almost local, $(\Gamma, \bar\sigma)$-invariant
elliptic pseudodifferential operator $Q$
and $(\Gamma, \bar\sigma)$-invariant smoothing operators
$R_0, \ R_1$ which satisfy
$$
QD = I - R_0 \qquad {\text{and}} \qquad DQ = I - R_1.
$$
\label{paramx}
\end{prop}

Define the idempotent
$$
   e(D) = \begin{pmatrix} {R}_0^2 & ({R}_0+
    {R}^2_0) {Q} \\ {R}_1{D} &
      1-{R}^2_1 \end{pmatrix}.
$$
Then $e(D)\in M_2({\mathcal R}(\Gamma, \sigma))$, where  ${\mathcal
R}(\Gamma, \sigma) =
\C(\Gamma, \sigma)\otimes\mathcal R$ is as
defined in \S 1.

The
${\mathcal R}(\Gamma, \sigma)$-{\em index} is by fiat
\[
   \Ind_\sigma(D)=[e(D)] - [E_0]\in
K_0({\mathcal R}(\Gamma, \sigma))\bigr)
\]
where $E_0$ is the idempotent
$$
E_0 = \begin{pmatrix} 0 & 0 \\ 0 &
      1 \end{pmatrix}.
$$
It is not difficult to see that $ \Ind_\sigma(D)$ is independent of the choice
of $(\Gamma, \bar\sigma)$-invariant parametrix $Q$ that is needed in its
definition.

Let $j: {\mathcal R}(\Gamma, \sigma) \to C^*_r(\Gamma, \sigma)$ be the
canonical
inclusion, which induces the morphism in $K$-theory
$$
j_* : K_\bullet({\mathcal R}(\Gamma, \sigma) ) \to
K_\bullet(C^*_r(\Gamma, \sigma)).
$$
Then we have,

\begin{defn*}
The $C^*_r(\Gamma, \sigma)$-{\em index} of a $(\Gamma, \bar\sigma)$-invariant
elliptic operator
$D : L^2(\M, \E) \to L^2(\M, \F)$ is defined as
$$
\Ind_{(\Gamma, \sigma)}(D) = j_*(\Ind_\sigma(D)) \in K_0(C^*_r(\Gamma, \sigma))
$$
\end{defn*}

\subsection{Heat kernels and the index map} Given $D$ as before, for
$t>0$,
we use the standard off-diagonal estimates for the heat kernel.
Recall that the heat kernels $e^{-tD^*D}$ and $e^{-tDD^*}$ are elements in the
${\mathcal R}(\Gamma, \sigma)$ (see the appendix).
Define the idempotent $e_t(D)\in M_2({\mathcal R}(\Gamma, \sigma))$
 (see the appendix) as follows
$$
   e_t(D) = \begin{pmatrix} e^{-t D^*D} & e^{-t/2 D^*D}\frac{(1-e^{-t D^*D})}
   {D^*D} D^* \\
e^{-t/2 D D^*}{D} &
      1- e^{-t DD^*} \end{pmatrix}.
$$
It is sometimes known as the Wasserman idempotent.

The relationship with the idempotent $e(D)$ constructed earlier can be
explained as follows. Define for $t>0$,
$$
Q_t = \frac{\left( 1-e^{-t/2 D^*D}\right)}{D^*D} D^*
$$
Then one easily verifies that $Q_t D = 1-e^{-t/2 D^*D} = 1- R_0(t)$
and $D Q_t  = 1-e^{-t/2 D D^*} = 1- R_1(t)$. That is, $Q_t$ is a parametrix
for $D$ for all $t>0$. Therefore one can write
$$
e_t(D) = \begin{pmatrix} R_0(t)^2 & (R_0(t) + R_0(t)^2) Q_t\\
R_1(t){D} &
     1- R_1(t)^2 \end{pmatrix}.
$$
In particular, one has for $t>0$
$$
\Ind_\sigma(D) = [e_t(D)] - [E_0] \in K_0({\mathcal R}(\Gamma, \sigma)).
$$

We use the same notation as in \cite{MM}.
A first order elliptic differential operator $D$ on $M$,
$$
D: L^2 (M, \mathcal E) \to L^2(M, \mathcal F)
$$
is by fiat a $\Gamma$-equivariant
first order
elliptic differential operator $\widetilde D$ on the smooth manifold
$\widetilde M$,
$$
\widetilde D :L^2 (\widetilde M, \E) \to L^2(\widetilde M, \F).
$$
Given any connection
$\nabla^{\W}$ on $\W$ which is compatible with the $\Gamma$ action and
the Hermitian metric, we define an extension of the
elliptic operator $\widetilde D$, to act
on sections of $\E\otimes \W$, $\F\otimes \W$,
\[
\widetilde D\otimes\nabla^{\W}:\Gamma(M,\E\otimes \W)\to
\Gamma(M,\F\otimes \W)
\]
as in \cite{MM}.

\subsection{Group cocycles and cyclic cocycles}

Using the pairing theory of cyclic cohomology and $K$-theory,
due to \cite{Co}, we will pair the $(\Gamma, \sigma)$-index of
a $(\Gamma, \bar\sigma)$-invariant elliptic operator $D$ on $\widetilde M$
with certain cyclic cocycles on ${\mathcal R}(\Gamma, \sigma)$. The
cyclic cocycles that we consider come from normalised group
cocycles on $\Gamma$. More precisely,
given a normalized group cocycle $c\in Z^k(\Gamma, \mathbb C)$, for $
k=0,\ldots, \dim M$,
we define a cyclic cocycle
$\tr_c$ of dimension $k$ on the twisted group ring
$\mathbb C(\Gamma, \sigma)$, which is
given by
$$
\tr_c(a_0\delta_{g_0}, \ldots,a_k\delta_{g_k}) = \left\{\begin{array}{l}
a_0\ldots a_k c(g_1,\ldots,g_k)
\tr( \delta_{g_0}\delta_{g_1} \ldots\delta_{g_k}) \quad\text{if} \,\, g_0
\ldots g_k =1\\
\\
0 \qquad \text{otherwise.}
\end{array}\right.
$$
where $a_j \in \C$ for $j=0,1,\ldots, k$.
To see that this is a cyclic cocycle on $\C(\Gamma, \sigma)$, we first
define,
as done in \cite{Ji}, the twisted differential graded algebra
$\Omega^\bullet(\Gamma, \sigma)$
as the differential graded algebra of finite linear combinations of
symbols
$$
g_0 dg_1\ldots dg_n \qquad g_i\in \Gamma
$$
with module structure and differential given by
\begin{align*}
(g_0 dg_1\ldots dg_n)g
&= \sum_{j=1}^n(-1)^{n-1} \sigma(g_j, g_{j+1}) g_0dg_1\ldots d(g_jg_{j+1})
\ldots dg_n  dg\\
& + (-1)^n \sigma(g_n, g) g_0 dg_1\ldots d(g_ng)\\
d(g_0 dg_1\ldots dg_n) & = dg_0 dg_1\ldots dg_n
\end{align*}
We now recall normalised group cocycles. A group $k$-cocycle is
a map
$h:\Gamma^{k+1} \to \C$ satisfying the identities
\begin{align*}
h(g g_0, \ldots g g_k) & = h(g_0, \ldots g_k)\\
0 & = \sum_{i=0}^{k+1} (-1)^i h(g_0, \ldots, g_{i-1},g_{i+1}\ldots, g_{k+1})
\end{align*}
Then a normalised group $k$-cocycle $c$ that is
associated to such an $h$ is given by
$$
c(g_1, \ldots, g_k) = h (1, g_1, g_1 g_2, \ldots, g_1\ldots g_k)
$$
and it is defined to be zero if either $g_i = 1$ or if $g_1\ldots g_k = 1$.
Any normalised group cocycle $c \in Z^k(\Gamma, \C)$
determines a $k$-dimensional
cycle via the following closed graded trace on
$\Omega^\bullet(\Gamma, \sigma)$
$$
\int g_0 dg_1\ldots dg_n =
\left\{\begin{array}{l}
 c(g_1,\ldots,g_k)
\tr( \delta_{g_0}\delta_{g_1} \ldots\delta_{g_k}) \quad\text{if} \, \, n=k
\quad\text{and}\,\, g_0
\ldots g_k =1\\
\\
0 \qquad \text{otherwise.}
\end{array}\right.
$$

Of particular interest is the case when $k=2$, when the formula above
reduces to
$$
\int g_0 dg_1 dg_2 =
\left\{\begin{array}{l}
 c(g_1, g_2)
\sigma(g_1, g_2) \quad \text{if} \,\, g_0 g_1 g_2 =1 ;\\
\\
0 \qquad \text{otherwise.}
\end{array}\right.
$$

The higher cyclic trace $\tr_c$ is by fiat this closed graded
trace.

\subsection{Twisted higher index theorem-
the cyclic cohomology version}

Let $M$ be a compact orbifold of dimension $n=4\ell$.
Let $\Gamma\to\widetilde{M}\overset{p}{\to}M$ be the universal
cover of $M$ and the orbifold fundamental group is $\Gamma$.
Let $D$ be an elliptic first order
operator on $M$ and $\widetilde{D}$ be the lift of
$D$ to $\widetilde{M}$,
\[
\widetilde{D}\,:\,L^2(\widetilde{M},
\E)\to
L^2(\widetilde{M},\F). \]
Note that $\widetilde{D}$ commutes with the $\Gamma$-action on
$\widetilde{M}$.

Now let $\omega$ be a closed 2-form on $M$ such that $\widetilde{\omega}
=p^* \omega=d\eta$ is \emph{exact} on $\widetilde M$. Define
$\nabla=d+\,i\eta$.
Then $\nabla$ is a Hermitian connection on the trivial line bundle
over $\widetilde{M}$, and the curvature of $\nabla,\ (\nabla)^2=i\,
\widetilde{\omega}$. Then $\nabla$ defines
a projective action of $\Gamma$ on $L^2$ spinors on $\widetilde M$ as in
section 1.

Consider the {\em twisted elliptic operator} on $\widetilde{M}$, \[
\widetilde{D}\otimes\nabla\,:\,L^2(\widetilde{M},
\E)\to
L^2(\widetilde{M},\F) \]
Then $\widetilde{D}\otimes\nabla$ no longer commutes with
$\Gamma$, but it does commute with the projective $(\Gamma,\bar\sigma)$ action.
In \S 2.1, we have defined the higher index of such an operator,
$$
\Ind_\sigma(\widetilde{D}\otimes\nabla)
\in K_0({\mathcal R}(\Gamma, \sigma)).
$$
Given a group cocycle $c\in Z^{2q}(\Gamma)$, one can define the
associated cyclic cocycle $\tau_c$ on ${\mathcal R}(\Gamma, \sigma)$
as in \S 2.3. Then $\tau_c$ induces a homomorphism on $K$-theory
$$
[\tau_c] : K_0({\mathcal R}(\Gamma, \sigma)) \to \R.
$$
The real valued higher index is the image of the higher index under this
homomorphism, i.e.
$$
\Ind_{(c, \Gamma, \sigma)}(\widetilde D\otimes \nabla)
= [\tau_c] (\Ind_\sigma(\widetilde{D}\otimes\nabla)  )
$$
To introduce the next theorem, we will briefly review some material on
characteristic classes for orbifold vector bundles. Let $M$ be a good
orbifold, that is the universal orbifold cover $\Gamma \to \widetilde M
\to M$ of $M$
is a smooth manifold. Then the orbifold tangent bundle $TM$ of $M$
can be
viewed as the $\Gamma$-equivariant bundle $T\widetilde M$ on $\widetilde M$.
Similar comments apply to the orbifold cotangent bundle $T^*M$ and, more
generally, to
any orbifold vector bundle on $M$. It is then clear that, choosing
$\Gamma$-invariant connections on the $\Gamma$-invariant vector bundles
on $\widetilde M$, one can define the Chern-Weil representatives of
the characteristic classes of the $\Gamma$-invariant vector bundles
on $\widetilde M$. These characteristic classes are $\Gamma$-invariant
and so define cohomology classes on $M$. For further details, see \cite{Kaw}.

\begin{thm}
Let $M$ be a compact, even dimensional, good
orbifold, and let $\Gamma$
be its orbifold fundamental group. Let $\widetilde D$ be a first order,
$\Gamma$-invariant elliptic differential operator acting on
$L^2$ sections of $\Gamma$-invariant vector bundles on $\widetilde M$, where
$\Gamma \to \widetilde M\to M$ is the universal orbifold cover of $M$.
Then, for any group cocycle $c\in Z^{2q}(\Gamma)$, one has
\begin{equation}
\Ind_{(c, \Gamma, \sigma)}(\widetilde D\otimes \nabla)
= \frac{q !} {(2\pi i)^q (2q!)} \left<Td(M)\cup ch(symb(D))
\cup f^*(\phi_c) \cup e^{\omega}, [T^*M]   \right>
\label{cohom} \end{equation}
where $Td(M)$ denotes the Todd characteristic class of
the complexified orbifold tangent bundle of $M$
which is pulled back to the orbifold cotangent bundle $T^*M$,
$ch(symb(D))$ is the Chern character of the symbol of the operator
$D$, $\phi_c$ is the Alexander-Spanier cocycle on $\underline B\Gamma$
that corresponds to the group cocycle $c$ and $f: M
\to \underline B\Gamma$ is the map that classifies the
orbifold universal cover  $ \M \to M$, cf. section 1.
\label{hindex}
\end{thm}

\begin{proof}
Choose a bounded, almost everywhere smooth Borel cross-section
$\beta : M\to \M$, which can then be used to define
the Alexander-Spanier cocycle $\phi_c$ corresponding to
$c\in Z^{2q}(\Gamma)$,  and such that $[\phi_c] = f^*[c]
\in H^{2q}(M)$. As in \S 2.2, for $t>0$, there is an
index idempotent,
$$
e_t(D) = \begin{pmatrix} R_0(t)^2 & (R_0(t) + R_0(t)^2) Q_t\\
R_1(t){D} &
     1- R_1(t)^2 \end{pmatrix}
\in M_2({\mathcal R}(\Gamma, \sigma)),
$$
where for $t>0$,
$$
Q_t = \frac{\left( 1-e^{-t/2 D^*D}\right)}{D^*D} D^*\quad, R_0(t) =
e^{-t/2 D^*D}, \quad R_1(t) = e^{-t/2 D D^*}.
$$
Then as in \S 2.2, one sees that $R_0(t)$, $R_1(t) $ are smoothing
operators and
$Q_t$ is a parametrix for $D$ for all $t>0$. The ${\mathcal R}(\Gamma,
\sigma)$-index map is then
\[
   \Ind_\sigma(\widetilde D\otimes \nabla)=
[e_t(\widetilde D\otimes \nabla)] - [E_0]\in
K_0({\mathcal R}(\Gamma, \sigma)).
\]
where $E_0$ is the idempotent
$$
E_0 = \begin{pmatrix} 0 & 0 \\ 0 &
      1 \end{pmatrix}.
$$
Let $R_t =e_t(\widetilde D\otimes \nabla) - E_0$. We adapt the strategy and
proof in \cite{CM} to our situation.
\begin{equation}\label{int}\begin{array}{lcl}
\Ind_{(c, \Gamma, \sigma)}(\widetilde D\otimes \nabla) &=&
\tr_c (R_t,R_t, \ldots R_t)\\[3mm]
&=& \displaystyle\int_{M^{2q+1}} \sum_{\gamma_1,\ldots\gamma_{2q}\in \Gamma}
c(1, \gamma_1,\ldots,\gamma_{2q})  \times \\[3mm]
& &{\tr}(R_t(\beta(x_0), \gamma_1\beta(x_1)) \ldots
R_t(\gamma_{2q}\beta(x_{2q}), \beta(x_0)))dx_0\ldots dx_{2q}
\end{array}\end{equation}
where $\beta: M \to \widetilde M$ denotes a bounded measurable section.
Notice that the phase term $\tr(\delta_{\gamma_1}\ldots
\delta_{\gamma_{2q}})$ appearing
in the expression for the cocycle $\tr_c$ cf. \S 2.3, is exactly cancelled
by the twisted
product of the integral kernels $R_t$, cf. \cite{CHMM} \S 3.
Notice also that the right hand side of equation (\ref{int}) is independent
of the
choice of the section $\beta$,  since upon changing $\beta$ to
$\gamma\cdot \beta$, we obtain {\small
$$
\begin{array}{lcl}
R_t(\gamma\beta(x_0), \gamma_1\gamma\beta(x_1)) \ldots
R_t(\gamma_{2q}\gamma\beta(x_{2q}), \gamma\beta(x_0)) &=&
e^{i\psi_\gamma(\beta(x_0))} R_t(\beta(x_0),
\gamma^{-1}\gamma_1\gamma
\beta(x_1))e^{-i\psi_\gamma(\gamma^{-1}\gamma_1\gamma \beta(x_1))}
\ldots \\[3mm]
& & e^{i\psi_\gamma(\gamma^{-1}\gamma_{2q}\gamma\beta(x_{2q}))}
R_t(\gamma^{-1}\gamma_{2q}\gamma\beta(x_{2q}),
\beta(x_0))e^{-i\psi_\gamma(\beta(x_0))} \\[3mm]
&=& R_t(\beta(x_0), \tilde\gamma_1\beta(x_1)) \ldots
R_t(\tilde\gamma_{2q}\beta(x_{2q}), \beta(x_0)),
\end{array}
$$}
where $\;\textsf{}\tilde\gamma_i=\gamma^{-1}\gamma_i \gamma\;$ which is
exactly as
in the case when the multiplier is trivial.

Observe that if $\beta_U : U \to \M$ is a smooth local
section, then there is a unique element $(\gamma_1, \ldots, \gamma_{2q})\in
\Gamma^{2q}$ such that
$(\beta(x_0), \gamma_1\beta(x_1),\ldots \gamma_{2q}\beta(x_{2q})) \in
\beta_U(U)^{2q+1}$. Moreover, we have the equality
$c(1, \gamma_1,\ldots,\gamma_{2q})  = \phi_c(x_0,x_1, \ldots x_{2q})$ \ \ (and
$\phi_c = 0$ otherwise),
where $\phi_c$ denotes the $\Gamma$-equivariant
(Alexander-Spanier) $2q$-cocycle on
$\widetilde M$ representing the pullback
$f^*(c)$ of the group 2-cocycle
via the classifying map $f$. Since
$R_t$ is mainly supported near the diagonal as $t\to 0$,
and using the equivariance of $R_t$, one sees that
$$
\Ind_{(c, \Gamma, \sigma)}(\widetilde D\otimes \nabla)
= \lim_{t\to 0}\int_{M^{2q+1}} \phi_c(x_0,x_1, \ldots x_{2q})
{\tr}(R_t(x_0, x_1)\ldots R_t(x_{2q}, x_0))dx_0dx_1\ldots dx_{2q},
$$
where we have identified $M$ with a fundamental domain for the
$\Gamma$ action on $\widetilde M$.
The proof is completed by applying the {\em local} higher index
Theorems 3.7 and 3.9 in \cite{CM}, to obtain the desired cohomological
formula (\ref{cohom}) for $\Ind_{(c, \Gamma, \sigma)}(\widetilde
D\otimes \nabla)$.
\end{proof}

\section{Twisted Kasparov map and range of the higher trace on $K$-theory}

In this section, we compute the range of the
2-trace $\tr_c$ on $K$-theory of the twisted group $C^*$ algebra,
where $c$ is a 2-cocycle on the group,
generalising the work of \cite{CHMM}.
Suppose as before that $\Gamma$ is a discrete, cocompact subgroup of
$PSL(2,\R)$
of signature $(g; \nu_1, \ldots, \nu_n)$.
That is, $\Gamma$ is the orbifold fundamental group of a
compact hyperbolic orbifold $\Sigma(g;\nu_1,\ldots,\nu_n)$
of signature $(g;\nu_1,\ldots,\nu_n)$.  Then for any multiplier
$\sigma$ on $\Gamma$ such that $\delta(\sigma) =0$, one has the
\emph{twisted Kasparov isomorphism},
\[
   \mu_\sigma : K^\bullet_{orb} (\Sigma(g;\nu_1,\ldots,\nu_n)) \to
K_\bullet (C^*_r(\Gamma,\sigma)),
\]
Proposition 2.14 in \cite{MM}. Its construction is recalled in this section,
as we need to refine it by factoring it through the $K$-theory
of the dense subalgebra  ${\mathcal R}(\Gamma,\sigma)$ of
$C^*_r(\Gamma,\sigma)$.
This is necessary in order to be able to use the pairing theory
of Connes \cite{Co}, \cite{CM} between higher cyclic traces
and $K$-theory.
We note that  using a
result of \cite{Ji}, that  ${\mathcal R}(\Gamma,\sigma)$ is
indeed a dense subalgebra of $C^*_r(\Gamma,\sigma)$ in our case.
In particular, given any
projection $P$ in $C^*_r(\Gamma,\sigma)$ there is
both a projection $\tilde P$ in the same $K_0$ class but lying in
the dense subalgebra ${\mathcal R}(\Gamma,\sigma)$.
This fact will also be utilized in the next section.
On the other hand, by the results of the current
section, given any such projection $P$ there
is a higher topological index that we can associate to it
cf. Theorem 3.3. The main result we prove
here is that the range of the 2-trace $\tr_c$ on $K$-theory of the
twisted group $C^*$ algebra is always an integer multiple
of a rational number. This will enable us to compute the range of values
of the Hall conductance in the quantum Hall effect on hyperbolic space,
generalizing the results in \cite{CHMM}.

\subsection{Twisted Kasparov map}
Let $\Gamma$ be as before,
that is, $\Gamma$ is the orbifold fundamental group of
the hyperbolic orbifold $\Sigma(g; \nu_1,\ldots,\nu_n)$.
Then for any multiplier $\sigma$ on $\Gamma$, we will factor
the \emph{twisted Kasparov isomorphism},
\begin{equation}\label{mu}
 \mu_\sigma : K^\bullet_{orb} (\Sigma(g; \nu_1,\ldots,\nu_n)) \to
K_\bullet (C^*_r(\Gamma,\sigma))
\end{equation}
in \cite{MM} through the $K$-theory
of the dense subalgebra  ${\mathcal R}(\Gamma,\sigma)$ of
$C^*_r(\Gamma,\sigma)$.

Let
$\mathcal{E}\to\Sigma(g; \nu_1,\ldots,\nu_n)$ be an orbifold
vector bundle over $\Sigma(g; \nu_1,\ldots,\nu_n)$
defining  an element $[\mathcal{E}]$
in $K^0(\Sigma(g; \nu_1,\ldots,\nu_n))$. As in \cite{Kaw}, one can form the
twisted Dirac
operator $\npartial^+_{\mathcal{E}} : L^2(\Sigma(g; \nu_1,\ldots,\nu_n),
\mathcal{S}^+\otimes\mathcal{E})
\to L^2(\Sigma(g; \nu_1,\ldots,\nu_n), \mathcal{S}^-\otimes\mathcal{E})$
where $\mathcal{S}^{\pm}$ denote the $\frac12$
spinor bundles over $\Sigma(g; \nu_1,\ldots,\nu_n)$.
One can lift the twisted Dirac operator $\npartial^+_{\mathcal{E}}$ as above,
to a $\Gamma$-invariant operator $\widetilde{\npartial_{\mathcal{E}}^+}$
on  $\mathbb{H} = \widetilde\Sigma(g; \nu_1,\ldots,\nu_n)$, which is
the universal orbifold cover of $\Sigma(g; \nu_1,\ldots,\nu_n)$,
\[
   \widetilde{\npartial_{\mathcal{E}}^+} : L^2(\mathbb{H},
      \widetilde{\mathcal{S}^+\otimes\mathcal{E}}) \to L^2(\mathbb{H},
      \widetilde{\mathcal{S}^-\otimes\mathcal{E}})
\]
For any multiplier $\sigma$ of $\Gamma$ with $\delta([\sigma]) = 0$,
there is a $\R$-valued 2-cocycle $\zeta$ on $\Gamma$ with
$[\zeta] \in H^2(\Gamma, \R)$ such that
$[e^{2\pi\sqrt{-1}\zeta}] = [\sigma]$. By the argument of \cite{MM},
section 2.2,  we know that we have an isomorphism
$H^2(\Gamma, \R) \cong H^2(\Gamma_{g'}, \R)$, and therefore
there is a 2-form $\omega$ on $\Sigma_{g'}$ such that
$[e^{2\pi\sqrt{-1}\omega}] = [\sigma]$.
Of course, the choice of $\omega$ is not unique, but this will not affect the
results that we are concerned with. Let $\widetilde \omega$ denote the lift
of $\omega$ to the universal cover $\mathbb H$. Since the hyperbolic plane
$\mathbb H$ is contractible, it follows that $\widetilde \omega = d\eta$
where $\eta$ is a 1-form on $\mathbb H$ which is not in general $\Gamma$
invariant. Now $\nabla = d + i\eta$ is a Hermitian connection on the trivial
complex line bundle on $\mathbb H$. Note that the curvature of $\nabla$ is
$\nabla^2 = i\tilde\omega$. Consider now the twisted Dirac operator
$\widetilde \npartial^+_{\mathcal{E}}$ which is
twisted again by the connection $\nabla$,
\[
   \widetilde{\npartial_{\mathcal{E}}^+}\otimes\nabla : L^2(\mathbb{H},
      \widetilde{\mathcal{S}^+\otimes\mathcal{E}}) \to L^2(\mathbb{H},
      \widetilde{\mathcal{S}^-\otimes\mathcal{E}}).
\]
It does not commute with the $\Gamma$ action, but it does commute
with the projective $(\Gamma, \bar\sigma)$-action  which is defined by the
connection
$\nabla$ as in \S 1. In section 2.1, we have defined the higher index of
such an
operator
\[
   \Ind_{\sigma} (\widetilde{\npartial^+_{\mathcal{E}}}\otimes
\nabla) \in
      K_0(\mathcal{R}(\Gamma, \sigma)),
\]
where  as before, $\mathcal{R}$ denotes the algebra of rapidly decreasing
sequences
on $\Z^2$.

Then the {\em twisted Kasparov map} (\ref{mu}) is
$$
\mu_\sigma([\mathcal{E}]) = j_*(\Ind_{\sigma}
(\widetilde{\npartial^+_{\mathcal{E}}}\otimes \nabla) )
= \Ind_{(\Gamma, \sigma)} (\widetilde{\npartial_{\mathcal{E}}^+}
\otimes\nabla)
\in
K_0({C}^*(\Gamma, \sigma)),
$$
where $j: {\mathcal R}(\Gamma, \sigma) = {\C}(\Gamma, \sigma) \otimes
\mathcal R  \to
{C}^*_r(\Gamma, \sigma)\otimes \mathcal K$
is the natural inclusion map, and as before, $\mathcal K$ denotes
the algebra of compact operators. Then
$$j_* : K_0 ( {\mathcal R}(\Gamma, \sigma)) \to K_0 ({C}^*_r(\Gamma,
\sigma))$$ is the
induced map on $K_0$. The twisted Kasparov map was defined for
certain torsionfree groups
in \cite{CHMM} and the general case in \cite{Ma1}.
It is related to the Baum-Connes assembly map \cite{BC}, \cite{BCH},
as is discussed in  \cite{Ma1}.

\subsection{Range of the higher trace on $K$-theory}
The first step in the proof is to show that
given a {\em bounded}
 group cocycle $c\in Z^2(\Gamma)$ we may define canonical pairings
with $K^0(\Sigma(g;\nu_1,\ldots,\nu_n))$ and $K_0(C^*_r(\Gamma, \sigma))$
which
are related by the twisted Kasparov isomorphism, by
adapting some of the results of Connes and
Connes-Moscovici to the twisted case.
As  $\Sigma(g;\nu_1,\ldots,\nu_n) = \underline B\Gamma$
is a negatively curved orbifold,
we know (by \cite{Mos} and \cite{Gr}) that degree 2 cohomology classes in
$H^2(\Gamma)$ have
{\em bounded} representatives i.e.\  bounded 2-cocycles on $\Gamma$.
The bounded group 2-cocycle $c$ may be
regarded as a skew symmetrised function on
$\Gamma\times\Gamma\times\Gamma$, so that
 we can use the results in section 2 to obtain a cyclic
2-cocycle $\tr_c$ on $\mathbb{C}(\Gamma, \sigma)\otimes \mathcal{R}$
by defining:
\[
   \tr_c(f^0\otimes r^0, f^1\otimes r^1, f^2\otimes r^2) =
{\mbox{Tr}}(r^0r^1r^2)
      \sum_{g_0g_1g_2=1} f^0(g_0) f^1(g_1) f^2(g_2) c(1, g_1, g_1g_2)
\sigma(g_1, g_2).
\]
Since the only difference
with the expression obtained in \cite{CM} is $\sigma(g_1, g_2)$, and since
$|\sigma(g_1, g_2)| = 1$, we can use Lemma 6.4, part (ii) in \cite{CM}
and the assumption that $c$ is bounded, to obtain the necessary
estimates which show that in fact $\tr_c$ extends continuously to the bigger
algebra ${\mathcal R}(\Gamma, \sigma)$. By
the pairing of cyclic theory and $K$-theory in \cite{Co},
one obtains an additive map
\[
   [\tr_c] : K_0(\mathcal{R}(\Gamma, \sigma))\to\mathbb{R}.
\]
Explicitly, $[\tr_c]([e]-[f]) = \widetilde\tr_c(e,\cdots,e) -
\widetilde\tr_c(f,\cdots,f)$, where $e,f$ are idempotent matrices with
entries in $(\mathcal{R}(\Gamma, \sigma))^\sim = $ unital
algebra obtained by adding the identity to
$\mathcal{R}(\Gamma, \sigma)$ and $\widetilde\tr_c$
denotes the
canonical extension of $\tr_c$ to $(\mathcal{R}(\Gamma, \sigma))^\sim$. Let \
$\widetilde{\npartial_{\mathcal{E}}^+}\otimes \nabla$ \  be
the Dirac operator defined in the previous section, which is
invariant under the projective action of the fundamental group
defined by $\sigma$.
Recall that by definition, the $(c,\Gamma, \sigma)$-index of
\ $\widetilde{\npartial_{\mathcal{E}}^+}\otimes \nabla$\ is
\[
   \Ind_{(c, \Gamma, \sigma)}(\widetilde{\npartial_{\mathcal{E}}^+}\otimes
\nabla)=
 [\tr_c](\Ind_{\sigma} (\widetilde{\npartial_{\mathcal{E}}^+}
\otimes\nabla)) = \langle[\tr_c], \mu_\sigma([\mathcal E])\rangle
\in\mathbb{R}.
\]
It only depends on the cohomology class $[c]\in H^2(\Gamma)$, and it is
linear with respect to $[c]$.  We assemble this to give the following theorem.

\begin{thm}
Given $[c] \in H^2(\Gamma)$ and $\sigma \in H^2(\Gamma, U(1))$ a
multiplier on $\Gamma$, there is a canonical additive map
\[
   \langle [c],\ \ \rangle :
K^0_{orb}(\Sigma(g;\nu_1,\ldots,\nu_n))\to\mathbb{R},
\]
which is defined as
\[
   \langle [c], [\mathcal{E}]\rangle =
\Ind_{(c, \Gamma, \sigma)}(\widetilde{\npartial_{\mathcal{E}}^+}\otimes
\nabla)=
 [\tr_c](\Ind_{\sigma} (\widetilde{\npartial_{\mathcal{E}}^+}
\otimes\nabla)) = \langle[\tr_c], \mu_\sigma([\mathcal E])\rangle
\in\mathbb{R}.
\]
Moreover, it is linear with respect to $[c]$.
\label{range}
\end{thm}

The {\em area cocycle} $c$ of the Fuchsian group $\Gamma$ is a canonically
defined 2-cocycle
on $\Gamma$ that is defined as follows.
Firstly, recall that
there is a well known area 2-cocycle on $PSL(2, \R)$, cf. \cite{Co2},
defined as follows:
$PSL(2, \R)$ acts on $\mathbb H$ such that $\mathbb H \cong PSL(2, \R)/SO(2)$.
Then $c(g_1, g_2) = \text{Area}(\Delta(o, g_1.o, {g_2}^{-1}.o)) \in \R$,
where $o$
denotes an origin in $\mathbb H$ and
$\text{Area}(\Delta(a,b,c))$ denotes the hyperbolic area of the
geodesic triangle in $\mathbb H$
with vertices at $a, b, c \in \mathbb H$. Then the restriction of
$c$ to the subgroup $\Gamma$ is the  area cocycle $c$ of $\Gamma$.

\begin{cor}
Let $c,\ [c]\in H^2(\Gamma)$, be the area cocycle, and ${\mathcal{E}}\to
\Sigma(g;\nu_1, \ldots, \nu_n)$ be an orbifold vector bundle over the orbifold
$\Sigma(g;\nu_1, \ldots, \nu_n)$.  Then in the notation above, one has
\[
   \langle [c], [{\mathcal{E}}]\rangle = \phi\rank\mathcal{E}
      \in\phi\mathbb{Z}.
\]
where $ - \phi = 2(1-g) + (\nu-n)\in {\mathbb Q}$ is the orbifold
Euler characteristic of $\Sigma(g;\nu_1, \ldots, \nu_n)$ and
$\quad \nu = \sum_{j=1}^n 1/\nu_j$.
\label{phi}
\end{cor}

\begin{proof}
By Theorem 2.2, one has
 one has
\begin{equation}
   [\tr_c](\Ind_{\sigma}
(\widetilde{\npartial_{\mathcal{E}}^+}\otimes
\nabla)) =
       \frac{1}{2\pi\#(G)}\int_{\Sigma_{g'}} \hat{A}(\Omega)
\tr(e^{R^{\mathcal{E}}}) e^{\omega}
\psi^*(\tilde c),\label{trc}
\end{equation}
where $\Sigma_{g'}$ is smooth and $G\to \Sigma_{g'} \to \Sigma(g;\nu_1,
\ldots, \nu_n)$
is a finite orbifold cover. Here
$\psi: \Sigma_{g'} \to \Sigma_{g'}$ is the lift of the
map $f:\Sigma(g; \nu_1, \ldots \nu_n) \to \Sigma(g; \nu_1, \ldots \nu_n)$
(since $\underline B\Gamma = \Sigma(g; \nu_1, \ldots \nu_n)$ in this case)
which is the classifying map of the orbifold universal
cover (and which in this case is the identity map) and $[\tilde c]$
degree 2 cohomology class on $\Sigma_{g'}$ that
is the lift of $c$ to $\Sigma_{g'}$.
We next simplify the right hand
side of (\ref{trc}) using the fact that $\hat{A}(\Omega) = 1$ and that
\begin{align*}
   \tr(e^{R^{\mathcal{E}}}) &= \rank \mathcal{E} + \tr({R^{\mathcal{E}}}), \\
   \psi^*(\tilde c) &= \tilde c,\\
   e^{\omega} &= 1 + {\omega} .
\end{align*}
We obtain
\[
   [\tr_c](\Ind_{\sigma}
(\widetilde{\npartial_{\mathcal{E}}^+}\otimes
\nabla)) =
      \frac{\rank \mathcal E}{2\pi\#(G)}
\langle [\tilde c], [\Sigma_{g'}] \rangle .
\]
When $c,\ [c]\in H^2(\Gamma)$, is the area 2-cocycle, then
$\tilde c$ is merely the restriction of the area cocycle on $PSL(2, \R)$
to the subgroup $\Gamma_{g'}$. Then one has
\[
   \langle [\tilde c], [\Sigma_{g'}] \rangle
= -2\pi\chi(\Sigma_{g'}) = 4\pi(g'-1).
\]
The corollary now follows from Theorem \ref{range} above together with
the fact that
$g' = 1 + \frac{\#(G)}{2}(2(g-1) + (n-\nu))$, and $\nu =
\sum_{j=1}^n 1/\nu_j$.
\end{proof}

We next describe the canonical pairing of $K_0(C^*_r(\Gamma, \sigma))$,
given $[c]\in H^2(\Gamma)$. Since $\Sigma(g;\nu_1,\ldots,\nu_n)$ is
negatively curved, we
know from \cite{Ji} that
\[
   {\mathcal R}(\Gamma, \sigma) = \left\{ f : \Gamma \to \mathbb{C} \mid
      \sum_{\gamma\in\Gamma} |f(\gamma)|^2 (1+l(\gamma))^k < \infty
      \mbox{ for all } k\ge 0\right\},
\]
where $l:\Gamma \to \mathbb{R}^+$ denotes the length function, is a dense
and spectral invariant subalgebra of $C_r^*(\Gamma, \sigma)$.  In
particular it is closed under the smooth functional calculus, and is
known as the algebra of rapidly decreasing $L^2$ functions on $\Gamma$.
By a theorem of \cite{Bost}, the inclusion map ${\mathcal R}(\Gamma,
\sigma)\subset
C^*_r(\Gamma, \sigma)$ induces an isomorphism
\begin{equation}
   K_j({\mathcal R}(\Gamma, \sigma)) \cong K_j(C_r^*(\Gamma, \sigma)),\quad
j=0,1.
\end{equation}
The desired pairing is the one obtained from the canonical pairing of
$K_0({\mathcal R}(\Gamma, \sigma))$ with $[c]\in H^2(\Gamma)$ using
the canonical isomorphism.  Therefore one has the equality
\[
   \langle [c], \mu_\sigma^{-1}[P]\rangle = \langle [\tr_c], [P] \rangle
\]
for any $[P]\in K_0({\mathcal R}(\Gamma, \sigma)) \cong K_0(C_r^* (\Gamma,
\sigma))$.
Using the previous corollary, one has

\begin{thm}[Range of the higher trace on $K$-theory]
Let $c$ be the area 2-cocycle on $\Gamma$.  Then $c$ is known to
be a bounded 2-cocycle, and one has
\[
   \langle [\tr_c], [P] \rangle = \phi(\rank\mathcal{E}^0 -
      \rank\mathcal{E}^1)\in\phi\mathbb{Z},
\]
where $ - \phi = 2(1-g) + (\nu-n)\in {\mathbb Q}$ is the orbifold Euler
characteristic of $\Sigma(g;\nu_1, \ldots, \nu_n)$
and $\quad \nu = \sum_{j=1}^n 1/\nu_j$.
Here $[P]\in K_0({\mathcal R}(\Gamma, \sigma)) \cong K_0(C_r^*(\Gamma,
\sigma))$,
and $\mathcal{E}^0, \ {\mathcal E}^1$ are orbifold vector bundles over
$\Sigma(g; \nu_1, \ldots,\nu_n)$ such that
\[
   \mu_\sigma^{-1}([P]) = [\mathcal{E}^0] -
      [{\mathcal E}^1] \in K^0_{orb}(\Sigma(g;\nu_1,\ldots,\nu_n)).
\]
In particular, the range of the the higher trace on $K$-theory is
$$
[\tr_c]( K_0 (C^*(\Gamma,\sigma))) = \phi\Z.
$$
\label{hrange}
\end{thm}

Note that
$\phi$ is in general only a {\em rational} number and we will give examples
to show that this is the case; however it is
an {\em integer} whenever the orbifold is smooth, i.e. whenever $1=\nu_1=\ldots
=\nu_n$, which is the case considered in \cite{CHMM}.
We will apply this result in the next section
to compute the range of values the Hall conductance in the quantum Hall effect
on the hyperbolic plane,
for orbifold fundamental groups, extending the results in \cite{CHMM}.

In the last section we provide a list of specific examples where
fractional values are achieved, and discuss the physical significance
of our model.

\section{The Area cocycle, the hyperbolic Connes-Kubo formula and
the Quantum Hall Effect}

In this section, we adapt and generalize the discrete model
of the quantum Hall effect of Bellissard
and his collaborators \cite{Bel+E+S} and also \cite{CHMM},
to the case of general cocompact Fuchsian
groups and orbifolds, which can  be viewed equivalently
as the generalization to the equivariant context.
We will first derive the discrete analogue
of the hyperbolic Connes-Kubo formula for the
Hall conductance 2-cocycle,
which was derived in the continuous case in \cite{CHMM}.
We then relate it
to the Area 2-cocycle on the twisted group algebra of the
discrete Fuchsian group, and we show that these define the same
cyclic cohomology class.
This enables us to use the results of the previous section to
show that the Hall conductance
has plateaux at all energy levels belonging to
any gap in the spectrum of the Hamiltonian, where it
is now shown to be equal to an integral multiple of
a {fractional} valued
topological invariant, namely the orbifold Euler characteristic.
The presence of denominators is
caused by the presence of cone points singularities
and by the hyperbolic
geometry on the complement of these cone points.
Moreover the set of possible denominators is finite and has been
explicitly determined in the next section, and the results
compared to the experimental data. It is plausible that this might
shed light on the mathematical mechanism responsible for
fractional quantum numbers in the quantum Hall effect.

We consider the Cayley graph
of the Fuchsian group $\Gamma$ of signature
$(g; \nu_1, \ldots, \nu_n)$, which acts freely on
the complement of a countable set of points in the
hyperbolic plane.
The Cayley graph embeds in the hyperbolic plane as follows.
Fix a base point $u \in \mathbb H$ such that the stabilizer (or isotropy
subgroup)
at $u$ is trivial and consider the orbit of the $\Gamma$ action through
$u$. This gives the vertices of the graph. The
edges of the graph are geodesics constructed as follows.
Each element of the group $\Gamma$
may be written as a word of minimal length in the
generators of $\Gamma$ and their inverses. Each generator
and its inverse determine a unique geodesic emanating from a vertex $x$
and these geodesics form the edges  of the
graph. Thus each word $x$ in the generators determines a
piecewise geodesic path from $u$ to $x$.

Recall that the {\em area cocycle} $c$ of the Fuchsian group
$\Gamma$ is a canonically defined 2-cocycle
on $\Gamma$ that is defined as follows.
Firstly, recall that
there is a well known area 2-cocycle on $PSL(2, \R)$, cf. \cite{Co2},
defined as follows:
$PSL(2, \R)$ acts on $\mathbb H$ such that
$\mathbb H \cong PSL(2, \R)/SO(2)$.
Then $c(\gamma_1, \gamma_2) = \text{Area}(\Delta(o, \gamma_1.o,
{\gamma_2}^{-1}.o)) \in \R$,
where $o$
denotes an origin in $\mathbb H$ and
$\text{Area}(\Delta(a,b,c))$ denotes the hyperbolic area of the
geodesic triangle in $\mathbb H$
with vertices at $a, b, c \in \mathbb H$. Then the restriction of
$c$ to the subgroup $\Gamma$ is the  area cocycle $c$ of $\Gamma$.

This area cocycle defines in a canonical way a cyclic 2-cocycle
$\tr_c$ on the group algebra $\C(\Gamma,\sigma)$ as follows;
$$
\tr_c(a_0,a_1,a_2) =
\sum_{\gamma_0\gamma_1\gamma_2 =1}
a_0(\gamma_0) a_1 (\gamma_1) a_2 (\gamma_2) c(\gamma_1, \gamma_2)
\sigma(\gamma_1, \gamma_2)
$$

We will now describe the hyperbolic Connes-Kubo formula for the Hall
conductance in the Quantum Hall Effect. Let $\Omega_j$ denote the (diagonal)
operator on $\ell^2(\Gamma)$ defined by
$$
\Omega_j f(\gamma) = \Omega_j(\gamma)f(\gamma) \quad \forall f\in
\ell^2(\Gamma)
\quad \forall \gamma\in \Gamma
$$
where
$$
\Omega_j(\gamma) = \int_o^{\gamma.o} \alpha_j \quad j=1,\ldots ,2g
$$
and where
\begin{equation} \label{sympl:basis}
\{\alpha_j \}_{j=1,\ldots ,2g}=\{ a_j \}_{j=1,\ldots ,g} \cup \{ b_{j}
\}_{j=1,\ldots ,g} \end{equation}
is a collection of harmonic $V$-forms on the orbifold
$\Sigma(g;\nu_1,\ldots,\nu_n)$, generating $H^1(\Sigma_g,\R)=\R^{2g}$,
cf. \cite{Kaw2} pg.78-83. These correspond to harmonic $G$-invariant
forms on $\Sigma_{g'}$ and to harmonic $\Gamma$-invariant forms on
$\mathbb H$.

Notice that we can write equivalently
$$
\Omega_j(\gamma) = c_j(\gamma),
$$
where the group cocycles $c_j$ form a symplectic basis for
$H^1(\Gamma,\Z)=\Z^{2g}$, with generators $\{ \alpha_j
\}_{j=1,\ldots ,2g}$, as in (\ref{sympl:basis}) and can be defined as
the integration on loops on
the Riemann surface of genus $g$ underlying the orbifold
$\Sigma(g;\nu_1,\ldots,\nu_n)$,
$$
c_j(\gamma)=\int_\gamma \alpha_j.
$$

For $j=1,\ldots ,2g$, define the derivations $\delta_j$ on
$\mathcal{R}(\Gamma,\sigma)$
as being the commutators $\delta_j a = [\Omega_j, a]$. A simple calculation
shows that
$$
\delta_j a (\gamma) = \Omega_j(\gamma)
a(\gamma) \quad \forall a\in \mathcal{R}(\Gamma,\sigma)
\quad \forall \gamma\in \Gamma.
$$

Thus, we can view this as the following general construction.
Given a 1-cocycle $a$ on the discrete group $\Gamma$,
i.e.
$$
a(\gamma_1\gamma_2) = a(\gamma_1) + a(\gamma_2) \qquad \forall
\gamma_1,\gamma_2\in \Gamma
$$
one can
define a derivation $\delta_a$ on the twisted group algebra
${\mathbb C}(\Gamma, \sigma)$
$$
\delta_a (f) (\gamma) = a(\gamma) f(\gamma).
$$
Then we verify that

\begin{align*}
\delta_a (fg) (\gamma) & = a(\gamma) fg(\gamma)\\
& = a(\gamma) \sum_{\gamma = \gamma_1\gamma_2} f(\gamma_1) g(\gamma_2)
\sigma(\gamma_1,\gamma_2)\\
& = \sum_{\gamma = \gamma_1\gamma_2} \Big(a(\gamma_1) + a(\gamma_2)\Big)
f(\gamma_1) g(\gamma_2) \sigma(\gamma_1,\gamma_2)\\
& = \sum_{\gamma = \gamma_1\gamma_2} \Big(\delta_a (f)(\gamma_1)
g(\gamma_2) \sigma(\gamma_1,\gamma_2)
+ f(\gamma_1) \delta_a (g) (\gamma_2) \sigma(\gamma_1,\gamma_2)\Big)\\
& = (\delta_a (f) g) (\gamma) + (f \delta_a g) (\gamma).
\end{align*}

As determined in section 1, the first cohomology of the group
$\Gamma=\Gamma(g;\nu_1,\ldots,\nu_n)$ is a free Abelian group of
rank $2g$. It is in fact a symplectic
vector space over $\mathbb Z$, and assume that $a_j, b_j, j=1,\ldots g$ is
a symplectic
basis of $H^1(\Gamma, \mathbb Z)$, as in (\ref{sympl:basis}). We denote
$\delta_{a_j}$ by $\delta_j$ and
 $\delta_{b_j}$ by $\delta_{j+g}$. Then these derivations give rise to
cyclic 2-cocycle
 on the twisted group algebra ${\mathbb C}(\Gamma, \sigma)$,
 $$
 \tr^K(f_0, f_1, f_2) = \sum_{j=1}^g \tr(f_0 (\delta_j(f_1)\delta_{j+g}(f_2) -
 \delta_{j+g}(f_1)\delta_j(f_2)))
 $$
$\tr^K$ is called the {\em Connes-Kubo Hall conductance cyclic
2-cocycle}.

In terms of the $\Omega_j$,
note that we have the simple estimate
$$
|\Omega_j(\gamma)|\le ||a_j||_{(\infty)} d(\gamma.o, o)
$$
where $d(\gamma.o, o)$ and the distance $d_\Gamma (\gamma, 1)$ in the
word metric on the group
$\Gamma$ are equivalent. This then yields the
estimate
$$
|\delta_j a (\gamma)| \le C_N d_\Gamma (\gamma, 1)^{-N}\quad \forall N\in
\mathbb N
$$
i.e  $\delta_j a \in \mathcal{R}(\Gamma,\sigma) \quad \forall a
\in \mathcal{R}(\Gamma,\sigma)$. Note that since $\forall \gamma,
\gamma' \in \Gamma$, the difference
$\Omega_j(\gamma\gamma') - \Omega_j(\gamma')$ is a constant
independent of $\gamma'$, we see that $\Gamma$-equivariance is preserved.
For $j=1,\ldots ,2g$, define the cyclic 2-cocycles
$$
\tr^K_j (a_0,a_1,a_2) = \tr(a_0 (\delta_ja_1 \delta_{j+g} a_2
- \delta_{j+g} a_1 \delta_j a_2)).
$$
These compute the Hall conductance for currents in the $(j+g)$th direction
which are induced by electric fields in the $j$th direction, as can be shown
using the quantum adiabatic theorem of Avron-Seiler-Yaffe \cite{Av+S+Y}
just as in section 6 of \cite{CHMM}, in the continuous model.
Then the {\em hyperbolic Connes-Kubo formula} for the Hall
conductance is the cyclic 2-cocycle given by the sum
$$
\tr^K (a_0,a_1,a_2) =
\sum_{j=1}^g \tr^K_j (a_0,a_1,a_2).
$$

\begin{thm}[The Comparison Theorem]
$$
[\tr^K] =  [\tr_c] \in HC^2(\mathcal{R}(\Gamma,\sigma))
$$
\label{comparison}
\end{thm}

\noindent{\bf Proof:}
Our aim is now to compare the two cyclic 2-cocycles
and to prove that they differ by a coboundary i.e.
$$
\tr^K (a_0,a_1,a_2) - \tr_c (a_0,a_1,a_2) =
b\lambda (a_0,a_1,a_2)
$$
for some cyclic 1-cochain $\lambda$ and where $b$ is the cyclic
coboundary operator. The key to this theorem is a
geometric interpretation
of the hyperbolic Connes-Kubo formula.

We begin with some calculations, to enable us to make this
comparison of the cyclic 2-cocycles.
$$
\tr^K (a_0,a_1,a_2) =
$$
$$
\sum_{j=1}^g \sum_{\gamma_0\gamma_1\gamma_2 =1}
a_0(\gamma_0) \left( \delta_j a_1(\gamma_1) \delta_{j+g} a_2(\gamma_2)
- \delta_{j+g} a_1(\gamma_1) \delta_j a_2(\gamma_2)\right)
\sigma(\gamma_0, \gamma_1) \sigma(\gamma_0\gamma_1, \gamma_2)
$$
$$
= \sum_{j=1}^g \sum_{\gamma_0\gamma_1\gamma_2 =1}
a_0(\gamma_0) a_1(\gamma_1) a_2(\gamma_2) \left(
\Omega_j(\gamma_1) \Omega_{j+g}(\gamma_2)
- \Omega_{j+g}(\gamma_1) \Omega_j(\gamma_2)
\right)\sigma(\gamma_1, \gamma_2)
$$
since by the cocycle identity for multipliers, one has
\begin{align*}
\sigma(\gamma_0, \gamma_1) \sigma(\gamma_0\gamma_1, \gamma_2)
& =  \sigma(\gamma_0,\gamma_1\gamma_2) \sigma(\gamma_1, \gamma_2)\\
& = \sigma(\gamma_0,\gamma_0^{-1}) \sigma(\gamma_1, \gamma_2) \quad
\text{since}
\quad \gamma_0\gamma_1\gamma_2 = 1\\
& = \sigma(\gamma_1, \gamma_2) \quad \text{since}
\quad\sigma(\gamma_0,\gamma_0^{-1}) =1.
\end{align*}
So we are now in a position to compare the two cyclic 2-cocycles.
Define $\Psi_j(\gamma_1, \gamma_2) =
\Omega_j(\gamma_1) \Omega_{j+g}(\gamma_2)
- \Omega_{j+g}(\gamma_1) \Omega_j(\gamma_2)$.

Let $\Xi: \mathbb H \to {\mathbb R}^{2g}$ denote the Abel-Jacobi
map
$$ \Xi: x \mapsto \left( \int_o^x a_1, \int_o^x
b_{1}, \ldots, \int_o^x a_g, \int_o^x b_{g} \right), $$
where $\displaystyle\int_o^x$ means integration along the unique geodesic in
${\mathbb H}$ connecting $o$ to $x$. The origin $o$ is chosen so that
it satisfies $\Gamma . o \cong\Gamma$.
The map $\Xi$ is a symplectic map, that is, if $\omega$ and $\omega_J$
are the respective symplectic
2-forms, then one has $\Xi^*(\omega_J) = \omega$.
One then has the following
geometric lemma.

\begin{lemma}
$$
\sum_{j=1}^g  \Psi_j(\gamma_1, \gamma_2) = \int_{\Delta_E(\gamma_1, \gamma_2)}
\omega_J
$$
where $\Delta_E(\gamma_1, \gamma_2)$ denotes the Euclidean triangle with
vertices at $\Xi(o), \Xi(\gamma_1.o)$ and $ \Xi(\gamma_2.o)$, and
$\omega_J$ denotes the flat K\"ahler 2-form on the Jacobi variety. That is,
$\sum_{j=1}^g  \Psi_j(\gamma_1, \gamma_2)$ is equal to the
Euclidean
area of the Euclidean triangle $\Delta_E(\gamma_1, \gamma_2)$.
\label{geom}
\end{lemma}

\begin{proof} We need to consider the expression
$$\sum_{j=1}^g  \Psi_j(\gamma_1, \gamma_2)
= \sum_{j=1}^g \Omega_j(\gamma_1) \Omega_{j+g}(\gamma_2)
- \Omega_{j+g}(\gamma_1) \Omega_j(\gamma_2).$$
Let $s$ denote the symplectic form on ${\mathbb R}^{2g}$ given by:
$$s(u,v)=\sum_{j=1}^g(u_jv_{j+g} -u_{j+g}v_j).$$
The so-called `symplectic area' of a triangle
with vertices $\Xi(o)=0,\Xi(\gamma_1.o),\Xi(\gamma_2.o)$ may be seen to be
$s(\Xi(\gamma_1.o),\Xi(\gamma_2.o))$. To appreciate this, however,
we need to use an argument from
\cite{GH}, pages 333-336.
In terms of the standard basis of ${\mathbb R}^{2g}$ (given in this case
by vertices in the integer period lattice arising from our choice of basis
of harmonic one forms) and corresponding coordinates
$u_1,u_2,\ldots u_{2g}$
the form $s$ is the two form on ${\mathbb R}^{2g}$ given by
$$\omega_J=\sum_{j=1}^g du_j\wedge du_{j+g}.$$
Now the  `symplectic area' of a triangle in ${\mathbb R}^{2g}$ with
 vertices $\Xi(o)=0,\Xi(\gamma_1.o),\Xi(\gamma_2.o)$
is given by integrating $\omega_J$ over the triangle and
a brief calculation reveals that this yields
$s(\Xi(\gamma_1.o),\Xi(\gamma_2.o))/2$, proving the lemma.
\end{proof}

We also observe that since $\omega = \Xi^* \omega_J$, one has
$$
c(\gamma_1, \gamma_2) = \int_{\Delta(\gamma_1, \gamma_2)}
\omega = \int_{\Xi(\Delta(\gamma_1, \gamma_2))}
\omega_J
$$
Therefore the difference
\begin{align*}
\sum_{j=1}^g  \Psi_j(\gamma_1, \gamma_2)  - c(\gamma_1, \gamma_2)
& = \int_{\Delta_E(\gamma_1, \gamma_2)}
\omega_J  - \int_{\Xi(\Delta(\gamma_1, \gamma_2))}
\omega_J\\
& = \int_{\partial \Delta_E(\gamma_1, \gamma_2)}
\Theta_J  - \int_{\partial\Xi(\Delta(\gamma_1, \gamma_2))}
\Theta_J
\end{align*}
where $\Theta_J$ is a 1-form on the universal cover of the
Jacobi variety such that $d\Theta_J = \omega_J$.  Therefore
one has
$$
\sum_{j=1}^g  \Psi_j(\gamma_1, \gamma_2)  - c(\gamma_1, \gamma_2)
= h(1, \gamma_1)  - h(\gamma_1^{-1}, \gamma_2)  +
h(\gamma_2^{-1}, 1)
$$
where $h(\gamma_1^{-1}, \gamma_2) =  \int_{\Xi(\ell (\gamma_1, \gamma_2))}
\Theta_J - \int_{m(\gamma_1, \gamma_2)}
\Theta_J $, where $\ell(\gamma_1, \gamma_2)$ denotes the
unique geodesic in $\mathbb H$ joining $\gamma_1.o$ and $\gamma_2.o$
and $m(\gamma_1, \gamma_2)$ is the straight line in the Jacobi variety
joining the points $\Xi(\gamma_1.o)$ and $\Xi(\gamma_2.o)$.
Since we can also write
$h(\gamma_1^{-1}, \gamma_2) = \int_{D(\gamma_1, \gamma_2)}
\omega_J
$, where $D(\gamma_1, \gamma_2)$ is a disk in the Jacobi variety
with boundary $\Xi(\ell(\gamma_1, \gamma_2))\cup m(\gamma_1, \gamma_2)$,
we see that $h$ is $\Gamma$-invariant.

We now define the cyclic 1-cochain $\lambda$ on $\mathcal{R}(\Gamma,\sigma)$
as
$$
\lambda(a_0, a_1) = \tr((a_{0})_h a_1) = \sum_{\gamma_0\gamma_1 =1}
h(1, \gamma_1) a_0(\gamma_0) a_1(\gamma_1)\sigma(\gamma_0, \sigma_1)
$$
where $(a_{0})_h$ is the operator on $\ell^2(\Gamma)$ whose matrix
in the canonical basis is
$h(\gamma_1, \gamma_2)a_0(\gamma_1\gamma_2^{-1})$.
Firstly, one has by definition
$$
b\lambda(a_0, a_1, a_2) =  \lambda(a_0 a_1, a_2)-
\lambda(a_0, a_1a_2)+\lambda(a_2a_0, a_1)
$$
We compute each of the terms seperately
\begin{align*}
\lambda(a_0 a_1, a_2) & =
\sum_{\gamma_0\gamma_1\gamma_2 =1}
h(1, \gamma_2) a_0(\gamma_0) a_1(\gamma_1)a_2(\gamma_2)
\sigma(\gamma_1, \gamma_2) \\
\lambda(a_0, a_1 a_2) & =
\sum_{\gamma_0\gamma_1\gamma_2 =1}
h(1, \gamma_1 \gamma_2) a_0(\gamma_0) a_1(\gamma_1)a_2(\gamma_2)
\sigma(\gamma_1, \gamma_2) \\
\lambda(a_2 a_0, a_1) & =
\sum_{\gamma_0\gamma_1\gamma_2 =1}
h(1, \gamma_1) a_0(\gamma_0) a_1(\gamma_1)a_2(\gamma_2)
\sigma(\gamma_1, \gamma_2)
\end{align*}
Now by $\Gamma$-equivariance, $h(1, \gamma_1\gamma_2)
= h(\gamma_1^{-1}, \gamma_2)$ and  $
h(1, \gamma_2) =  h(\gamma_2^{-1}, 1)$. Therefore one has
$$
b\lambda(a_0, a_1, a_2) =
$$
$$
\sum_{\gamma_0\gamma_1\gamma_2 =1}
 a_0(\gamma_0) a_1(\gamma_1)a_2(\gamma_2)
\left(h(\gamma_2^{-1}, 1) - h(\gamma_1^{-1}, \gamma_2)
+ h(1, \gamma_1)\right)
\sigma(\gamma_1, \gamma_2)
$$
Using the formula above, we see that
$$
b\lambda(a_0, a_1, a_2) =   \tr^K (a_0,a_1,a_2) - \tr_c (a_0,a_1,a_2).
$$

It follows from Connes pairing theory of cyclic cohomology and $K$-theory
\cite{Co2}, by the range of the higher trace Theorem \ref{hrange} and
by the Comparison Theorem 4.1 above that

\begin{cor}[Rationality of conductance]
The Connes-Kubo Hall conductance cocycle $\tr^K$ is rational.
More precisely, one has
$$
\tr^K (P, P, P) = \tr_c (P, P, P) \in \phi\mathbb Z
$$
for all projections $P\in \mathcal{R}(\Gamma,\sigma)$, where
$ -\phi = 2(1-g) + (\nu-n)\in {\mathbb Q}$ is the orbifold Euler
characteristic of $\Sigma(g;\nu_1, \ldots, \nu_n)$.
\label{trKtrc}
\end{cor}

Finally, suppose that we are given a very thin sample
of pure metal, with electrons situated along
the Cayley graph of $\Gamma$, and a very strong
magnetic field which is uniform and normal in direction to the sample.
Then at very low temperatures, close to absolute zero,
quantum mechanics dominates and the discrete model
that is considered
here is a model of electrons moving on the
Cayley graph of $\Gamma$ which is
embedded in the sample.
The associated discrete Hamiltonian $H_\sigma$
for the electron in the magnetic field
is given by the Random Walk operator in the
projective $(\Gamma, \sigma)$ regular representation on the
Cayley graph of the group $\Gamma$. It is also known as the generalized
{\em Harper operator} and was first studied in this generalized
context in \cite{Sun}, see also \cite{CHMM}. We will see that
the Hamiltonian that we consider is in a natural way
the sum of
a free Hamiltonian and a term that models the Coulomb interaction.
We also add a restricted class of potential terms to the
Hamiltonian in our model.

Because the charge carriers are Fermions, two different charge carriers
must occupy different quantum eigenstates of the Hamiltonian.
In the limit of zero temperature they minimize the energy and occupy
eigenstates with energy lower that a given one, called the
{\em Fermi level}
and denoted $E$. Let $P_E$ denote denote the
corresponding  spectral projection of the Hamiltonian.
If $E$ is not in the spectrum of
the Hamiltonian, then then $P_E \in
\mathcal{R}(\Gamma,\sigma)$ and
the hyperbolic Connes-Kubo formula for the
Hall conductance $\sigma_E$ at the energy level $E$ is
defined as follows;
$$
\sigma_E = \tr^K (P_E, P_E, P_E).
$$
As mentioned earlier, it measures the sum of the contributions
to the Hall conductance at the energy level $E$
for currents in the $(j+g)$th direction which are induced by
electric fields in the $j$th direction, cf. section 6 \cite{CHMM}.
By Corollary 4.3,
one knows that the Hall conductance takes on values
in $\phi \mathbb Z$  whenever
the energy level $E$ lies in a gap in the spectrum of the Hamiltonian
$H_\sigma$. In fact we notice that the Hall conductance is a constant function
of the energy level $E$ for all values of $E$ in the same gap
in the spectrum of the Hamiltonian. That is, the Hall conductance
has plateaux which are {\em integer} multiples of the fraction $\phi$
on the gap in the spectrum of the Hamiltonian.

We now give some details. Recall the left $\sigma$-regular representation
$$
(U(\gamma) f)(\gamma') = f(\gamma^{-1}\gamma')
\sigma(\gamma', \gamma^{-1}\gamma')
$$
$\forall f\in \ell^2(\Gamma)$ and $\forall \gamma, \gamma' \in \Gamma$.
It has the property that
$$
U(\gamma) U(\gamma') = \sigma(\gamma, \gamma')
U(\gamma\gamma')
$$
Let $S= \{A_j, B_j, A_j^{-1}, B_j^{-1}, C_i, C_i^{-1} : j= 1, \ldots, g, \quad
i=1, \ldots,n\}$ be a
symmetric set of generators for $\Gamma$. Then the Hamiltonian
is explicitly given as
\begin {align*}
H_\sigma & :\ell^2(\Gamma) \to \ell^2(\Gamma)\\
H_\sigma & = \sum_{\gamma\in S} U(\gamma)
\end{align*}
and is clearly by definition a bounded self adjoint operator.
Notice that the Hamiltonian can be decomposed as a sum of a free
Hamiltonian containing the torsionfree generators
and a term simulating Coulomb interactions, that contains the
torsion generators.
$$
H_\sigma  = H_\sigma^{free} + H_\sigma^{interaction}
$$
where
$$
H_\sigma^{free} = \sum_{j=0}^g U(A_j) +  U(B_j) + \left(U(A_j) +
U(B_j)\right)^*
$$ and
$$
H_\sigma^{interaction}  =\sum_{i=1}^n U(C_i) + U(C_i)^*.
$$
Let $V \in \C(\Gamma, \sigma)$ be any "potential", and
$$H_{\sigma, V} = H_\sigma + V.$$

\begin{lemma}
If $E\not\in \text{spec}(H_{\sigma, V} )$, then $P_E \in
\mathcal{R}(\Gamma,\sigma)$,
where $P_E = \chi_{[0,E]}(H_{\sigma, V} )$ is the spectral
projection of the Hamiltonian to energy levels less than or
equal to $E$.
\label{specproj}
\end{lemma}

\begin{proof}
Since $E\not\in \text{spec}(H_{\sigma, V} )$, then
$P_E = \chi_{[0,E]}(H_{\sigma, V} ) = \varphi(H_{\sigma, V} )$ for some smooth,
compactly supported function $\varphi$. Now by definition, $H_\sigma \in
\C(\Gamma,\sigma) \subset  \mathcal{R}(\Gamma,\sigma)$,
and since $\mathcal{R}(\Gamma,\sigma)$ is closed under the smooth
functional calculus by the result of \cite{Ji}, it follows that
$P_E \in \mathcal{R}(\Gamma,\sigma)$.

\end{proof}

Therefore by Corollary 4.3
and the discussion following it, we have,

\begin{thm}[Fractional Quantum Hall Effect] Suppose that the Fermi
energy level
$E$ lies in a gap of the spectrum
of the Hamiltonian $H_{\sigma, V} $, then the Hall conductance
$$
\sigma_E = \tr^K (P_E, P_E, P_E) = \tr_c (P_E, P_E, P_E)  \in \phi \mathbb Z
$$
That is, the Hall conductance
has plateaux which are {\em integer} multiples of $\phi$
on any gap in the spectrum of the Hamiltonian, where
$ -\phi = 2(1-g) + (\nu-n)\in {\mathbb Q}$ is the orbifold Euler
characteristic of $\Sigma(g;\nu_1, \ldots, \nu_n)$.
\label{gqhe}
\end{thm}

\begin{rems}
The set of possible
denominators $\phi$ for low genus coverings
can be derived easily from the results of \cite{Bro} and
is reproduced in the second table in the next section.
It is plausible that this Theorem
might shed light on the mathematical mechanism
responsible for fractional quantum numbers that occur
in the Quantum Hall Effect, as we attempt to explain in the
following
section.

\end{rems}

\section{Fractional quantum numbers: phenomenology}

We first discuss the characteristics of our model explaining the
appearance of fractional quantum numbers in the quantum Hall
effect. In particular, we point out the main advantages and
limitations of the model.

Our model is a single electron model. It is well known that the FQHE
is a consequence of the Coulomb interaction between electrons, hence
it should not be seen by a single particle model. However, in our
setting, the negative curvature of the hyperbolic structure provides a
geometric replacement for interaction. The equivalence between
negative curvature and interaction is well known from the case of
classical mechanics where the Jacobi equation for a single particle
moving on a negatively curved manifold can be interpreted as the
Newton equation for a particle moving in the presence of a negative
potential energy \cite{Arn}.

The main advantage of this setting is that the fractions derived in this
way are topological. In fact, they are obtained from an equivariant
index theorem. Moreover, they are completely determined by the geometry of the
orbifold. In fact, we have
$$ \phi=-\chi_{orb}(\Sigma(g;\nu_1,\ldots,\nu_n)). $$
Let us recall that the {\em orbifold Euler characteristic}
$\chi_{orb}(\Sigma)$ of an orbifold $\Sigma$, is a rational
valued invariant that is completely
specified by the following properties, cf. \cite{Tan}:
\begin{enumerate}
\item it is multiplicative under
orbifold covers;
\item it coincides with the topological Euler
characteristic in the case of a smooth surface;
\item it satisfies the
volume formula,
$$ \chi_{orb}(\Sigma_1 \cup\cdots \cup\Sigma_k)=\sum_{j=1}^k
\chi_{orb}(\Sigma_j) -\sum_{i,j} \chi_{orb}(\Sigma_i\cap
\Sigma_j)+\cdots (-1)^{k+1}\chi_{orb}(\Sigma_1\cap\cdots\cap
\Sigma_k),$$
whenever all the intersections on the right hand side are suborbifolds
of $\quad\Sigma_1 \cup\cdots \cup\Sigma_k$, and all the $\Sigma_j$ are
orbifolds of the same dimension.
\end{enumerate}
This characterization allows for ease of computation and prediction of
expected fractions.

Most notably, as pointed out in \cite{Bel+E+S}, the topological nature of
the Hall conductance makes it stable under small deformations of
the Hamiltonian. Thus, this model can be easily generalized to
systems with disorder, cf. \cite{CHM}. This is a necessary step in
order to establish the presence of plateaux \cite{Bel+E+S}.

The identifications of fractions with integer multiples of the orbifold
Euler characteristic imposes some restrictions on the range of
possible fractions from the geometry of the orbifolds. For instance,
it is known from the Hurwitz theorem that the maximal order of a
finite group acting by isometries on a smooth Riemann surface
$\Sigma_{g'}$ is $\#(G)=84(g'-1)$. Moreover, this maximal order is
always attained. Thus, the smallest possible fraction that appears in
our model is $\phi= \frac{2(g'-1)}{84(g'-1)}=1/42$.

This is, in some respects, an advantage of the model, in as it gives very
clear prediction on which fractions can occur, and at the same time
its main limitation, in as we do not get a complete agreement between
the set of fractions we obtain and the fractions that are actually
observed in experiments on the FQHE.

In order to compare our predictions with experimental data, we
restrict our attention to orbifolds with a torus or a sphere as
underlying topological surface.
Recall that, as explained above, we think of the hyperbolic
structure induced by the presence of cone points on these surfaces as a
geometric way of introducting interaction in this single electron
model, hence we would consider equivalently the underlying surface
with many interacting electrons (fractions observed in FQHE
experiments) or as a hyperbolic surface with one electron.

We report a table of comparison between the values obtained
experimentally and our prediction. Notice how the fraction $5/2$
which appears in the experimental values and caused major problems of
interpretation in the many-particle models appears here naturally as
the orbifold Euler characteristic of $\Sigma(1;6,6,6)$ (which we may
as well refer to as the Devil's orbifold).

{\small

\begin{tabular}{|c||c|} \hline
experimental &   $g =1$ or  $g=0$  \\ \hline \hline
$5/3$ & $\Sigma(1;6,6)$     \\  \hline
$4/3$ & $\Sigma(1;3,3)$    \\ \hline
$7/5$ & $\Sigma(0;5,5,10,10)$  \\ \hline
$4/5$ & $\Sigma(1;5)$   \\ \hline
$5/7$ & $\Sigma(0;7,14,14)$ \\ \hline
$2/3$ & $\Sigma(1;3)$   \\ \hline
$3/5$ & $\Sigma(0;5,10,10)$   \\ \hline
$4/7$ & $\Sigma(0;7,7,7)$ \\ \hline
$5/9$ &  ??? \\ \hline
$4/9$ & $\Sigma(0;3,9,9)$ \\ \hline
$3/7$ &  ???  \\ \hline
$2/5$ & $\Sigma(0;5,5,5)$ \\ \hline
$1/3$ & $\Sigma(0;3,6,6)$ \\ \hline
$5/2$ & $\Sigma(1;6,6,6)$  \\ \hline
\end{tabular}

}

Despite the small number of discrepancies in the table above,
the agreement between values of orbifold Euler characteristics and
experimentally observed fractions in the quantum Hall effect
is far from being satisfactory. In particular, not only there is a
small number of observed values which are not orbifold Euler
characteristics, but there are also many rational numbers that are
realized as orbifold Euler characteristics, which do not seem to
appear among the experimental data. For instance, by looking at the
values of the next table, reported also in figure 1, we see clearly
that we have some fractions with even denominator, such as $1/4$,
$1/2$, and $1/6$, which do not correspond to experimental values.
As pointed out in the
introduction, the reason for this discrepancy is that a more
sophisticated model for the Coulomb interaction is needed in general.

In the remaining
of this section, we discuss some phenomenology,
with particular emphasis on the nature of the
cone points and the role of
the minimal genus of the covering surface $\Sigma_{g'}$. We hope to
return to these topics in some future work.

Every orbifold $\Sigma(g;\nu_1,\ldots,\nu_n)$
is obtained as a quotient of a surface $\Sigma_{g'}$ with respect to
the action of a finite group $G$, cf. \cite{Sc}. In general both $g'$ and
$G$ are not unique. For instance, the orbifold $\Sigma(1;2,2)$ is
obtained as the quotient of $\Sigma_2$ by the action of $\Z_2$,
or as the quotient of $\Sigma_3$ by the action of $\Z_4$, or by the
action of $\Z_2\times \Z_2$, cf. \cite{Bro}.
For every $\Sigma(g;\nu_1,\ldots,\nu_n)$ there is a minimal
$g'$ such that the orbifold is obtained as a quotient of $\Sigma_{g'}$
by a finite group action. In \cite{Bro},
Broughton has derived a complete list of all the good two
dimensional orbifolds which are quotients of Riemann surfaces
$\Sigma_{g'}$ with genus  $g' = 2$ or $3$.

In a physical model one
can distinguish between two types of {\em disorder}: a mobility
disorder and a sample disorder,  cf. \cite{Bel+E+S}.
We can argue phenomenologically that, if an orbifold can
be realized by a covering of low genus, this corresponds to a lower
density of atoms in the sample, as opposed to the case of a surface of
high genus, as one can see by looking at the Cayley graph of
$\Gamma_{g'}$. Thus, we can consider the minimal genus of the smooth
coverings as a measure of mobility.
This means that, in an experiment, the fractions derived from
orbifolds with low genus coverings will be easier to observe (have
more clearly marked plateaux) than fractions which are only realized
by quotients of surfaces of higher genus.

Thus, we can consider the list of examples given in \cite{Bro} and
compute the corresponding fractions.
We list the result in the following table.

{\small

\begin{tabular}{|c|c||c|} \hline
$\phi$ &   $g' =2$ &   $g'=3$  \\ \hline \hline
$4/3$  &   & $\Sigma(0;3,3,3,3,3)$  $\Sigma(1;3,3)$ \\  \hline
$2/3$ &  $\Sigma(0;3,3,3,3)$ & $\Sigma(0;2,2,6,6)$
$\Sigma(0;2,3,3,6)$ $\Sigma(0;2,2,2,2,3)$  $\Sigma(1;3)$ \\ \hline
$4/7$ & &  $\Sigma(0;7,7,7)$ \\ \hline
$1/2$ & $\Sigma(0;2,2,4,4)^*$  $\Sigma(0;2,2,2,2,2)^*$ &
$\Sigma(0;4,8,8)$
$\Sigma(1;2)$ \\ \hline
$4/9$ & & $\Sigma(0;3,9,9)$ \\ \hline
$2/5$ & $\Sigma(0;5,5,5)$ & \\ \hline
$1/3$ & $\Sigma(0;3,6,6)$ $\Sigma(0;2,2,3,3)^*$ &
$\Sigma(0;2,12,12)$ $\Sigma(0;3,4,12)$
$\Sigma(0;4,4,6)$ $\Sigma(0;2,2,2,6)$ \\ \hline
$1/4$ & $\Sigma(0;2,8,8)^*$ $\Sigma(0;4,4,4)^*$
$\Sigma(0;2,2,2,4)^*$ &  \\ \hline
$1/5$ & $\Sigma(0;2,5,10)$ & \\ \hline
$4/21$ & & $\Sigma(0;3,7,7)$ \\ \hline
$1/6$ & $\Sigma(0;3,4,4)^*$ $\Sigma(0;2,6,6)^*$
$\Sigma(0;2,2,2,3)^*$ & $\Sigma(0;2,4,12)$ $\Sigma(0;3,3,6)$ \\ \hline
$1/8$ & $\Sigma(0;2,4,8)^*$ & \\ \hline
$1/12$ & $\Sigma(0;2,4,6)^*$ $\Sigma(0;3,3,4)^*$ & \\ \hline
$1/24$ & $\Sigma(0;2,3,8)^*$ & \\ \hline
$1/42$ & & $\Sigma(0;2,3,7)$ \\ \hline
\end{tabular}

}

\begin{figure}[ht]
\epsfig{file=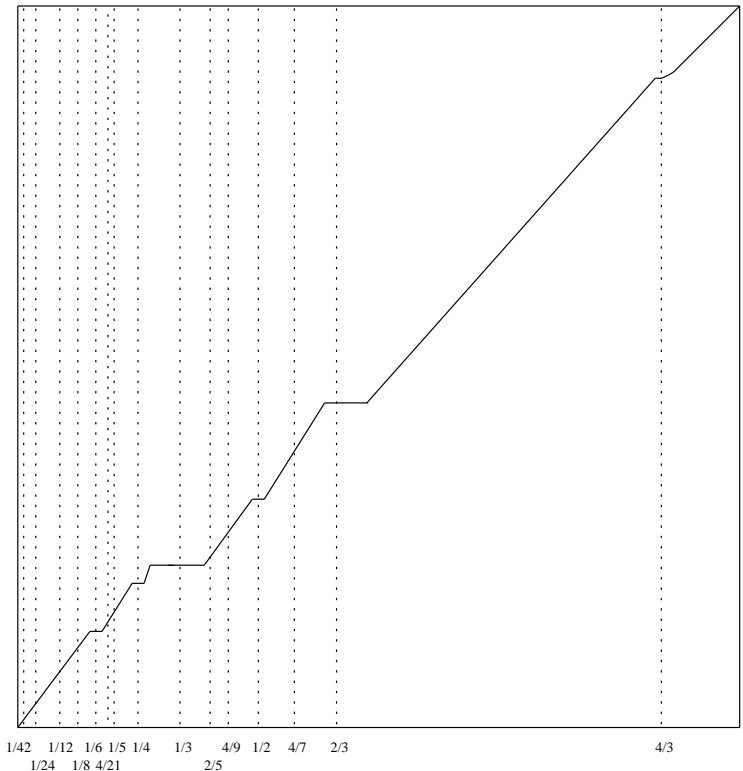,angle=270}
\caption{Phenomenology of fractions in the low genus approximation
\label{fig1} }
\end{figure}

In the table the orbifolds that are markes with a $*$ can
be realized both as quotient of $\Sigma_2$ and of $\Sigma_3$.
It seems also reasonable to think that if the same fraction is realized by
several different orbifolds, for fixed $g'$, then the corresponding
plateau will be more clearly marked in the experiment. This
would make $\phi=1/3$ the most clearly pronounced plateau, which is in
agreement with the experimental data. However, higher genus
corrections are not always negligible. In fact, by only considering
genus $g'=2$ and $g'=3$ contributions, we would expect a more marked
plateau for the fraction $\phi=2/3$ than for the fraction $\phi=2/5$,
and the experimental results show that this is not the case.
It seems important to observe that this model produces equally easily
examples of fractions with odd or even denominators (e.g. $\phi=1/4$
appears in the table above). It is
interesting to compare this datum with the difficulty encountered
within other models in explaining the appearance of the
fraction $5/2$ in the experiments. Its presence is only
justified by introducing a different physical model (the so called
non-abelian statistics).
In figure \ref{fig1} we sketch the plateaux as they would appear in the
result of an experiment, using only the low genus $g'=2$ and $g'=3$
approximation.

As we already mentioned in the introduction, both the
hyperbolic structure and  the
cone points are essential in order to have fractional quantum
numbers. In fact, $\phi$ is an
{\em integer} whenever the hyperbolic orbifold is smooth, i.e. whenever
$1=\nu_1=\ldots =\nu_n$, which is the case considered in
\cite{CHMM}. On the other hand, by direct inspection, it is possible
to see that all euclidean orbifolds also produce only {\em integer}
values of $\phi$. (Notice that sometimes hyperbolic orbifolds with cone
points may still produce integers: the orbifold $\Sigma(1;2,2)$ has
$\phi=1$, cf. \cite{Bro}.) Models of FQHE on euclidean orbifolds
have been considered, in a different, string-theoretic context,
e.g. \cite{Sk-Th}.

We can argue that the cone points can also be thought of as a form of
``disorder''. In fact, we may identify the preimage of the cone points
in the universal covering ${\mathbb H}$ with sample disorder (with
respect to
the points in the Cayley graph of $\Gamma_{g'}$). The same
fraction can often be obtained by orbifolds with a varying number of
cone points (for fixed $g'$), as illustrated in the previous
table. This can be
rephrased by saying that the system allows for more or less sample
disorder, and in some cases this can be achieved without affecting the
mobility measured by $g'$.

\section*{Appendix}

The main purpose of this appendix is to establish Lemma F, which is used
in the paper. We follow closely the approach in \cite{BrSu}.
We use  the notation of the previous sections.
Let $A$ be an operator on $L^2(\widetilde M, \widetilde S\otimes E)$
with Schwartz kernel  $k_A$ and also
commuting with the given
$(\Gamma ,{\bar\sigma})$-action.
Then one has
\begin{equation}
\label{A1}
e^{i\phi_\gamma(x)} k_A(\gamma x,\gamma y)\,
e^{-i\phi_\gamma(y)}
= k_A(x,y)\qquad\forall\gamma\in\Gamma,
\end{equation}
where we have identified the fibre at $x\in \widetilde M$
with the fibre at $\gamma x\in \widetilde M$.
If $k_A$ is smooth, then one can define the {\em von Neumann trace}
just as Atiyah did in the untwisted case,
\[
\tr\left(k_A\right) =
\int_{\mathcal F}\,\underline\tr \left(k_A(x,x)\right)\,dx, \]
where $\mathcal F$ denotes a fundamental domain for the action
of $\Gamma$ on $\widetilde M$ and
where $\underline\tr$ denotes
the pointwise or local trace. The von Neumann trace
is well defined, since as a consequence of (\ref{A1}),
$\;\underline\tr(k_A(x,x))$ is a $\Gamma$-invariant function on
$\widetilde{M}$. The following lemma establishes that it is a
trace.

\begin{lemma*}[A]
Let $A, B$ be operators on $L^2(\widetilde M, \widetilde S\otimes E)$
with smooth Schwartz kernels and also
commuting with the given
$(\Gamma,{\bar\sigma})$-action.
Then one has
$$
\tr\left(A B\right) = \tr\left(B A\right).
$$
\end{lemma*}

\begin{proof}
Let $k_A, k_B$ denote the smooth Schwartz kernels of $A,B$ respectively, and
$k_{AB}, k_{BA}$ denote the smooth Schwartz kernels of $AB,BA$
respectively. Then one has
\begin{align*}
\tr\left(A B - B A\right) & =
\int_{x\in \mathcal F} \underline\tr\left(k_{AB}(x,y) -
k_{BA}(x,y)\right)\\
& = \int_{x\in \mathcal F} \int_{y\in \widetilde M}\underline\tr\left(
k_{A}(x,y)k_{B}(y,x) -  k_{B}(x,y)k_{A}(y,x)\right)\\
& = \sum_{\gamma\in \Gamma}\int_{x\in \mathcal F} \int_{y\in \mathcal F}
\underline\tr\left( k_{A}(x,\gamma y)k_{B}(\gamma y,x) -
k_{B}(x,\gamma y)k_{A}(\gamma y,x)\right)\\
& = 0
\end{align*}
since each term in the summand vanishes by symmetry, and we have used
the fact that the fundamental domain $\mathcal F$ is compact in order
to interchange the order of the summation and integral.
\end{proof}

We will also adopt a more operator theoretic approach.
Let $\mathcal H =
L^2(\mathcal F, \widetilde S\otimes E|_{\mathcal F})$. Then
$\Phi: L^2(\widetilde M, \widetilde S\otimes E)
\stackrel{\cong}{\to} \ell^2(\Gamma, \mathcal H)$ is given by
$(\Phi s)(\gamma) = R_{\mathcal F}(T_\gamma s) \; \forall
\gamma \in \Gamma$, where $R_{\mathcal F}:
L^2(\widetilde M, \widetilde S\otimes E) \to \mathcal H$
denotes the restriction map to the fundamental domain $\mathcal F$.

As in section 1, let $W^*(\sigma) $ denote the commutant, i.e.
$$
W^*(\sigma) = \left\{
A \in \ell^2(\Gamma, \mathcal H) : [T_\gamma, A] = 0 \quad
\forall \gamma\in \Gamma\right\}
$$
Then one has the following simple lemma,

\begin{lemma*}[B]
$W^*(\sigma)$ is a semifinite von Neumann algebra.
\end{lemma*}

\begin{proof}
We need to show that $W^*(\sigma) $ is a $*$-algebra
which is weakly closed. We will establish that it is has a semifinite
trace a bit later on.

Let $A, B \in W^*(\sigma) $. Since $[T_\gamma, AB]
= [T_\gamma, A]B + A[T_\gamma, B]$, it follows that $AB \in
W^*(\sigma) $.
Since $[T_\gamma, A] = -[T_\gamma^*, A^*] = -[T_{\gamma^{-1}}, A^*]$
it follows $A^* \in W^*(\sigma) $.
Clearly the identity operator is in $W^*(\sigma) $.
Finally, if $A_n \in W^*(\sigma)  \; \forall n\in \mathbb N$
and $A_n$ converges weakly to $A$, it follows that
for all $\gamma\in \Gamma$, $T_\gamma A_n$
converges weakly to $T_\gamma A$ and also to $A T_\gamma$. By
uniqueness of weak limits, we deduce that $A \in
W^*(\sigma) $.
\end{proof}

For $A \in W^*(\sigma) $, define its generalized {\em Fourier
coefficients} $\widehat A(\gamma)\in B(\mathcal H)$ as
$$
\widehat A(\gamma) v = T_\gamma(A\delta_1^v)(1)
$$
where $\delta_1^v \in \ell^2(\Gamma, \mathcal H)$ is defined
for all $v\in\mathcal H$ as
$$
\delta_1^v (\gamma) = \left\{\begin{array}{l}
 v \;\;  {\rm if} \;\; \gamma=1;\\[7pt]
 0 \;\;  {\rm otherwise.}
\end{array}\right.
$$
Since  $\;T_\gamma\delta_1^v (\gamma') = \delta_1^v
(\gamma'\gamma)\sigma(\gamma',\gamma)$, one has
$$
T_\gamma\delta_1^v (\gamma') = \left\{\begin{array}{l}
 v \;\;  {\rm if} \;\; \gamma'=\gamma^{-1};\\[7pt]
 0 \;\;  {\rm otherwise,}
\end{array}\right.
$$
since $\sigma(\gamma^{-1}, \gamma) =1\; \forall \gamma\in\Gamma$.
In particular, it follows that for all $f \in \ell^2(\Gamma, \mathcal H)$,
one has
$$
f(\gamma) = \sum_{\gamma_1 \gamma_2 = \gamma} T_{\gamma_1}
\delta_1^{f(\gamma_2)}
$$
so that one has the following {\em Fourier expansion}
$$
\begin{array}{lcl}
 A f (\gamma) & = & \displaystyle\sum_{\gamma_1 \gamma_2 = \gamma}
AT_{\gamma_1} \delta_1^{f(\gamma_2)}
 = \displaystyle\sum_{\gamma_1 \gamma_2 = \gamma} T_{\gamma_1} A
\delta_1^{f(\gamma_2)} \\[+7pt]
&= &\displaystyle\sum_{\gamma_1 \gamma_2 =
\gamma}\widehat A(\gamma_1)(f(\gamma_2)).
\end{array}
$$
The following
elementary properties are satisfied by the Fourier coefficients.

\begin{lemma*}[C] For $A, B \in W^*(\sigma) $ and for all
$\gamma\in \Gamma$, for all $f \in \ell^2(\Gamma, \mathcal H)$, one has
\begin{itemize}
\item[(1)] ${A} f (\gamma) =
\displaystyle\sum_{\gamma_1 \gamma_2 = \gamma}\widehat A(\gamma_1)(
f(\gamma_2))$;
\item[(2)]${\widehat A^*}(\gamma) = (\widehat A(\gamma^{-1}) )^*$;
\item[(3)] ${\widehat {A B}} (\gamma) =
\displaystyle\sum_{\gamma_1 \gamma_2 = \gamma}\widehat A(\gamma_1)\widehat
B(\gamma_2)$;
\item[(4)] ${\widehat {A A^*}} (1) = \displaystyle\sum_{\gamma} \widehat
A(\gamma)
\widehat A(\gamma)$;
\item[(5)] $||A|| \le \displaystyle\sum_{\gamma}
||\widehat A(\gamma)||$;
\item[(6)] ${\widehat{A - B}}(\gamma) =
\widehat{A}(\gamma) - \widehat{B}(\gamma)$.
\end{itemize}
\end{lemma*}

\begin{proof} The proof follows by straightforward calculations as done
above. The reader is warned that the righthand side of the inequality in
part $(5)$ is not necessarily finite.
\end{proof}

Define $C_o(\Gamma, \mathcal K)$ to be the set of all
$A \in W^*(\sigma) $
such that $\widehat A(\gamma) \in \mathcal K\;\;\forall \gamma\in\Gamma$, and
$\widehat A(\gamma) = 0$ for all but finitely many $\gamma\in\Gamma$. Then
the completion of $C_o(\Gamma, \mathcal K)$ with respect to the operator norm
is denoted, as in section 1 of \cite{MM}, by $C^*_r(\Gamma, \sigma)
\otimes\mathcal K
$, and
called the twisted crossed product algebra associated to the twisted action
$(\alpha, \sigma)$. Then one has the following useful containment criterion,

\begin{lemma*}[D]
If $A \in W^*(\sigma)$ and also satisfies
$
\;\displaystyle\sum_{\gamma} ||\widehat A(\gamma)||<\infty,
$ then $A \in C^*_r(\Gamma, \sigma)
\otimes\mathcal K $.

If $A \in W^*(\sigma)$ and also satisfies
$
\;\displaystyle\sum_{\gamma} d(\gamma, 1)^k
||\widehat A(\gamma)||<\infty,
$ for all positive integers $k$, then $A \in  {\mathcal R}(\Gamma, \sigma)
$.

\end{lemma*}

\begin{proof}
Let $K_1 \subset K_2 \subset \cdots $ be a sequence of finite subsets of
$\Gamma$
which is an exhaustion of $\Gamma$, i.e. $\bigcup_{j\ge 1} K_j = \Gamma$.
For all
$j \in \mathbb N$, define $A_j \in W^*(\sigma)$ by
$$
{\widehat {A_j}} (\gamma) = \left\{\begin{array}{l}
{\widehat {A}} (\gamma) \;\;  {\rm if} \;\; \gamma\in K_j;\\[7pt]
 0 \;\;  {\rm otherwise.}
\end{array}\right.
$$
Then in fact $A_j \in C_o(\Gamma, \mathcal K)$ by definition, and
using the previous lemma, we have
$$
\begin{array}{lcl}
||A-A_j|| & \le & \displaystyle\sum_{\gamma}
||{\widehat{A - A_j}}(\gamma)||\\[+7pt]
& = & \displaystyle\sum_{\gamma}
||\widehat{A}(\gamma) - \widehat{A_j}(\gamma)||\\[+7pt]
& = & \displaystyle\sum_{\gamma\in \Gamma\setminus K_j}
||\widehat A(\gamma)||.
\end{array}
$$ By hypothesis, $\displaystyle\sum_{\gamma}
||\widehat A(\gamma)||<\infty$, therefore
$\displaystyle\sum_{\gamma\in \Gamma\setminus K_j}
||\widehat A(\gamma)|| \to 0$
as $j\to\infty$, since $K_j$ is an increasing exhaustion
of $\Gamma$. This
proves that
$A \in  C^*_r(\Gamma, \sigma)
\otimes\mathcal K $.

The second part is clear from the definition, once we
identify $\mathcal R$ with
the algebra of sequences
$$
\left\{\left(a_\gamma\right)_{ \gamma\in \Gamma}\Big|
\sup_{\gamma\in
\Gamma} d(\gamma, 1)^k |a_\gamma| < \infty\; \forall k\in
\mathbb N\right\}.
$$
\end{proof}

The following off-diagonal estimate is well known, cf. \cite{BrSu}.

\begin{lemma*}[E]
Let $D= \widetilde{\not\partial}^+_E\otimes\nabla^s$ be
a twisted Dirac operator. Then the Schwartz kernel $k(t,x,y)$ of the
heat operator $e^{-tD^*D}$ is smooth $\;\;\forall
t>0$. It also satisfies the following off-diagonal estimate
$$
|k(t,x,y)| \le C_1 t^{-n/2} e^{-C_2 d(x,y)^2/t}
$$ uniformly in $(0,T] \times\widetilde M\times\widetilde M$ for any $T>0$,
where $d$ denotes the Riemannian distance function on $\widetilde M$.
The same result is true for the the Schwartz kernel of the
heat operator $e^{-tDD^*}$.
\end{lemma*}

\begin{lemma*}[F]
Let $D= \widetilde{\not\partial}^+_E\otimes\nabla^s$ be
a twisted Dirac operator.
Then $e^{-tD^*D}, e^{-tDD^*} \in
{\mathcal R}(\Gamma, \sigma) \subset C^*_r(\Gamma, \sigma)
\otimes\mathcal K \;\;\forall
t>0$.
\end{lemma*}

\begin{proof}
By the Lemma above, it follows that $e^{-tD^*D}, e^{-tDD^*}$ are bounded
operators
commuting with the given twisted action, i.e. $e^{-tD^*D},\; e^{-tDD^*} \in
W^*(\Gamma, \mathcal H)$. Since the Schwartz kernels of
$\widehat{e^{-tD^*D}}(\gamma),\; \widehat{e^{-tDD^*}}(\gamma)$ are
smooth $\forall \gamma\in\Gamma$ by the Lemma above, it follows that
$\widehat{e^{-tD^*D}}(\gamma), \widehat{e^{-tDD^*}}(\gamma) \in \mathcal K$
$\forall \gamma\in\Gamma$. Let $d_\Gamma$ denote the word metric with
respect to a given finite set of generators, and $d$ the Riemannian metric on
$\widetilde M$. Then it is well known that
$$
d_\Gamma(\gamma_1, \gamma_2) \le C_3 (\inf_{x,y\in\widetilde M} d (\gamma_1
x, \gamma_2 y)
+1)
$$ for some positive constant $C_3$. By the Lemma 5 above, one has,
$$
||\widehat{e^{-tD^*D}}(\gamma)|| \le C_4 e^{-C_5 d_\Gamma(\gamma, 1)^2}
$$ for some positive constants $C_4, C_5$, and a similar estimate holds
for $\widehat{e^{-tDD^*}}(\gamma)$.
Setting $r(\gamma) = d_\Gamma(\gamma, 1)$ observe that one has the estimate
$$
\# \left\{\gamma\in\Gamma\; |\; r(\gamma) \le R\right\} \le C_6 e^{C_7 R}
$$
for some positive constants $C_6, C_7$, since the volume
growth rate of $\Gamma$ is at most exponential.
Therefore one has
$$
\sum_\gamma d(\gamma, 1)^k ||\widehat{e^{-tD^*D}}(\gamma)|| <\infty
\qquad{\rm and}
\qquad\sum_\gamma  d(\gamma, 1)^k
||\widehat{e^{-tDD^8}}(\gamma)|| <\infty
$$
for all positive integers $k$.
By the Lemma above, it follows that
$\;e^{-tD^*D}, e^{-tDD^*} \in
{\mathcal R}(\Gamma, \sigma) \subset C^*_r(\Gamma, \sigma)
\otimes\mathcal K \;\;\forall
t>0$.

\end{proof}


\begin{thebibliography}{Av+K+P+S3s}

\bibitem[Arn]{Arn} V.I. Arnold, Mathematical methods of
classical mechanics, {\em Graduate Texts in Mathematics}, Vol.{\bf 60},
Springer 1978.

\bibitem[At]{At} M.F. Atiyah,
Elliptic operators, discrete groups and Von Neumann algebras, {\em
Ast\'erisque}
{\bf 32-33} (1976), 43-72.

\bibitem[Av+S+Y ]{Av+S+Y}{J. Avron, R. Seiler, I. Yaffe},
{Adiabatic theorems and applications to the integer quantum Hall effect},
{\em Commun. Math. Phys.} {\bf 110} (1987), {33-49}.


\bibitem[BC]{BC} P. Baum and A. Connes,
Chern Character for discrete groups, in
A fete of Topology, Academic Press (1988) 163-232.

\bibitem[BCH]{BCH} P. Baum, A. Connes and N. Higson,
Classifying space for proper actions and $K$-theory of group
$C^*$-algebras,
{\em Contemp. Math.} {\bf 167} (1994) 241-291.

\bibitem[Bel+E+S]{Bel+E+S}{J. Bellissard, A. van Elst, H. Schulz-Baldes},
{The non-commutative geometry of the quantum Hall effect},
{\em J. Math. Phys.} {\bf 35} (1994), {5373-5451}.


\bibitem[Bost]{Bost} J. Bost, Principe d'Oka, $K$-th\'eorie et syst\'emes
dynamiques non
commutatifs, {\em Invent. Math.} {\bf 101} (1990), no. 2, 261-333.

\bibitem[Bro]{Bro} A. Broughton, Classifying finite group actions
on surfaces of low genus, {\em J. Pure Appl.
Algebra} {\bf 69} (1991), no. 3, 233-270.

\bibitem[BrSu]{BrSu} J. Br\"uning, T. Sunada, On the spectrum of
gauge-periodic elliptic
operators. M\'ethodes semi-classiques, Vol. 2 (Nantes, 1991). {\em
Ast\'erisque} {\bf 210} (1992),
65-74.

\bibitem[CHMM]{CHMM} A. Carey, K. Hannabuss, V. Mathai and P. McCann,
Quantum Hall Effect on the hyperbolic plane, {\em Commun. Math.
Physics}, {\bf 190} no. 3 (1998) 629-673.

\bibitem[CHM]{CHM} A. Carey, K. Hannabuss, V. Mathai, Quantum Hall
effect on the hyperbolic plane in the presence of disorder,
{\em Lett. Math. Phys.} {\bf 47} (1999), no. 3, 215--236


\bibitem[Co]{Co} A. Connes, Non commutative differential
geometry, {\em Publ. Math. I.H.E.S.} {\bf 62} (1986), 257-360.

\bibitem[Co2]{Co2} A. Connes, {\em Noncommutative geometry},
Academic Press, Inc., San Diego,
CA, (1994).

\bibitem[CM]{CM} A. Connes, H. Moscovici, Cyclic cohomology, the Novikov
conjecture
and hyperbolic groups, {\em Topology} {\bf 29} (1990), 345-388.

\bibitem[Far]{Far} C. Farsi, $K$-theoretical index theorems for good orbifolds,
{\em Proc. Amer. Math. Soc.} {\bf 115} (1992) 769-773.

\bibitem[Froh]{Froh} J. Frohlich, Transport in thermal equilibrium,
gapless modes, and anomalies, in
Festschrift for the 40th anniversary of the IHES,
{\em Publ. Math. I.H.E.S.} (1998)  81-97.

\bibitem[FuSt]{FuSt} M. Furuta and B. Steer, Seifert fibred homology
$3$-spheres and the Yang-Mills equations
on Riemann surfaces with marked points, {\em Adv. Math.} {\bf 96} (1992),
no. 1, 38-102.

\bibitem[GH]{GH} P. Griffiths and J. Harris, {\em Principles of
algebraic geometry}, Wiley, New York, 1978.

\bibitem[Gr]{Gr} M. Gromov, Volume and bounded cohomology,
{\em Publ. Math. I.H.E.S.} {\bf 56} (1982),  5-99.

\bibitem[Gr2]{Gr2} M. Gromov, K\"ahler-hyperbolicity and $L^2$ Hodge theory,
{\em J. Diff. Geom.}  {\bf  33}   (1991),   263-292.

\bibitem[Ji]{Ji} R. Ji, Smooth dense subalgebras of reduced group
$C^*$-algebras, Schwartz
cohomology of groups and cyclic cohomology, {\em Jour. Func. Anal.} {\bf
107} (1992), 1-33.

\bibitem[Kaw]{Kaw} T. Kawasaki, The index of elliptic operators over
$V$-manifolds,
{\em Nagoya Math. Jour.} {\bf 84} (1981) 135-157.

\bibitem[Kaw2]{Kaw2} T. Kawasaki, The signature theorem for
$V$-manifolds, {\em Topology}, {\bf 17} (1978), no. 1, 75--83.

\bibitem[MM]{MM} M. Marcolli, V. Mathai, Twisted index theory on good
orbifolds, I:
noncommutative Bloch theory, {\em Communications in
Contemporary Mathematics} {\bf 1} no. 4 (1999) 553-587.




\bibitem[Ma1]{Ma1} V. Mathai, K-theory of twisted group $C^*$-algebras and
positive scalar curvature,
Rothenberg Festschrift, {\em Contemp. Math.} {\bf 231} (1999)
203-225.


\bibitem[Mos]{Mos} G. Mostow, {\em Strong rigidity of symmetric
spaces}, Ann. Math. Studies,
{\bf 78} (1973), Princeton University Press.

\bibitem[Patt]{Patt} S.J. Patterson, On the cohomology of Fuchsian
groups, {\em Glasgow Math. J.} {\bf 16} (1975), no. 2, 123--140.


\bibitem[Sc]{Sc} P. Scott, The geometries of 3-manifolds, {\em Bull. Lond.
Math. Soc.} {\bf 15} (1983) 401-487.

\bibitem[Si]{Si} I.M. Singer, Some remarks on operator theory
and index theory, in
$K$-theory and operator algebras,
Lecture Notes in Math.,
Vol. 575, Springer, Berlin, (1977)  128-138.

\bibitem[Sk-Th]{Sk-Th} S. Skoulakis, S. Thomas, Orbifold duality
symmetries and quantum Hall systems, {\em Nucl.Phys. B} {\bf 538}
(1999) 659-684.


\bibitem[Sun]{Sun} T. Sunada,
{A discrete analogue of periodic magnetic Schr\"odinger operators},
{\em Contemp. Math.} {\bf 173} (1994), 283-299.

\bibitem[Tan]{Tan} C. Tanasi,
The Euler--Poincar\'e characteristic of two-dimensional orbifolds,
{\em Rend. Sem. Mat. Univ. Politec. Torino} {\bf 45} (1987) 133-155.



\bibitem[Xia]{Xia} J. Xia, Geometric invariants of the quantum Hall effect,
{\em Commun. Math. Phys.} {\bf 119} (1988), 29-50.

\end{thebibliography}
\end{document}